\documentclass[a4paper,11pt]{article}
\usepackage{graphicx,graphics,xcolor}
\usepackage{epsfig}
\usepackage{amsmath}
\usepackage{amssymb}
\usepackage{mathtools}
\usepackage[top=3cm, bottom=3.5cm, left=2cm, right=2cm]{geometry}
\usepackage{microtype}
\usepackage{lmodern}
\usepackage[utf8]{inputenc}
\usepackage{bera}
\newcommand{\ds}{\displaystyle }

\newcommand{\beq}{\begin{equation} }
\newcommand{\eeq}{\end{equation}}
\author{{\Large Thomas M. Michelitsch$^a$, Federico Polito$^b$,
Alejandro P. Riascos$^c$ }\\ \\ 
\footnotesize{$^a$ Sorbonne Universit\'e, Institut Jean le Rond d’Alembert, 
CNRS UMR 7190} \\
\footnotesize{4 place Jussieu, 75252 Paris cedex 05, France} \\ 
\footnotesize{E-mail: michel@lmm.jussieu.fr}\\[1ex]
\footnotesize{$^b$ Department of Mathematics ``Giuseppe Peano'', University of Torino, Italy}  \\ 
\footnotesize{E-mail: federico.polito@unito.it}\\[1ex]
\footnotesize{$^c$Instituto de F\'isica, Universidad Nacional Aut\'onoma de M\'exico} \\[1ex]
\footnotesize{Apartado Postal 20-364, 01000 Ciudad de M\'exico, M\'exico}  \\
\footnotesize{E-mail: aperezr@fisica.unam.mx} \\[4ex]
}
\title{Squirrels can remember little: A random walk with jump reversals induced by a discrete-time renewal process}
\begin{document}

\maketitle
\begin{abstract}
\footnotesize{
We consider a class of discrete-time random walks with directed unit steps on the integer line.
The direction of the steps is reversed at the time instants of events in a  discrete-time renewal process and is maintained at uneventful time instants. 
This model represents a discrete-time semi-Markovian generalization of the telegraph process.
We derive exact formulae for the propagator using generating functions. We 
prove that for geometrically distributed waiting times in the diffusive limit, this walk converges to the classical telegraph process.
We consider the large-time asymptotics of the expected position: For waiting time densities with finite mean the walker remains in the average localized close to the departure site whereas escapes for fat-tailed waiting-time densities (i.e. densities with infinite mean) by a sublinear power-law. 
We explore anomalous diffusion features by accounting for the `aging effect' as a hallmark of non-Markovianity
where the discrete-time version of the `aging renewal process' comes into play. By deriving pertinent distributions
of this process we obtain explicit formulae for the variance when the waiting-times are Sibuya-distributed. In this case and generally for fat-tailed waiting time PDFs a $t^2$-ballistic superdiffusive scaling emerges in the large time limit.
In contrast if the waiting time PDF between the step reversals is light-tailed (`narrow' with finite mean and variance) the walk exhibits normal diffusion and for `broad' waiting time PDFs (with finite mean and infinite variance) superdiffusive large time scaling.
We also consider time-changed versions where the walk is subordinated to a continuous-time point process such as the time-fractional Poisson process. This defines a new class of biased continuous-time random walks exhibiting several regimes of anomalous diffusion. 
}
 \\[1mm]
\end{abstract}
\newpage
\tableofcontents

\section{Introduction}
\label{Introduction}
The topic of biased random walks on the integer line has a long history 
\cite{Doob1953,Feller1971,SpitzerF1976,Hajek1987}. Originally the main motivation stems from the gambler's ruin problem and related contexts such as betting. In the last few decades this interest is enhanced by the upswing of approaches involving fractional calculus with applications in anomalous transport and diffusion  \cite{MetzlerKlafter2001,GorenfloMainardi2008,Gorenflo2009,MainardiRobertoScalas2000,
ScalasGorenfloMainardi2004,MainardiGorScalas2004,Laskin2003,GorenfloMainardi2011,BeghinOrsingher2009,Meerschaert2011,fractional_book2019} and generalized fractional dynamics \cite{CahoyPolito2013,MichelitschRiascos2020a,MichelitschRiascos2020b,Sandev2018} (and many others). Most of these models are based on the continuous-time random walk (CTRW) approach of Montroll and Weiss \cite{MontrollWeiss1965}, including
asymmetric anomalous transport \cite{ADTRW2021,WangBarkai2020} 
and dynamics in networks \cite{RiascosMateos2021,Barabasi2016,Newman2018} (and many others).
A special kind of 
biased walk with a complete memory of its history is the so called ``elephant random
walk'' (ERW) introduced in 2004 by Sch\"utz and Trimper \cite{SchuetzTrimper2004} and consult \cite{BaurBertoin2020} for the remarkable connection with P\'olya urns. 
\\[1mm]
The present paper is devoted to studying a random walk on the one-dimensional infinite lattice where directed unit steps are performed at integer time instants. The walker changes the step direction at the event instants of a discrete-time renewal process and maintains the step direction at uneventful time instants.
Comparing this definition with the ERW \cite{SchuetzTrimper2004}, our
walker remembers only the last decision taken at the latest event of the renewal process. Therefore, in a sense our walker has a weaker memory as the `elephant' walker.
Invoking the popular belief that squirrels sometimes forget where they buried their nuts (having a weaker memory as elephants which ``always remember''), we shall refer to our random walker as `squirrel' and call our walk `{\it squirrel random walk}' (SRW).
\\[1mm]
A descriptive picture for the ERW are indeed P\'olya urns (see \cite{BaurBertoin2020} for details). 
The  P\'olya urn picture also shows the main difference between the elephant and the squirrel walker. The elephant walker
considers at each integer time instant the entire history of the previous steps where the numbers
of steps done in one direction (red urns) and in the opposite direction (white urns) govern the
probability for the random choice of the new step direction. Contrarily to the elephant the squirrel
walker does not ‘remember’ the entire set of previous steps. Compared to elephants, squirrels
can remember little, namely only the choice of the step direction at the latest renewal time instant in
a discrete-time counting process. The weaker memory of squirrels compared to elephants
therefore does not allow to establish a simple connection with P\'olya urns as an appropriate picture for the SRW.
\\[1mm]
Indeed, it turns out that the SRW is a semi-Markovian
discrete-time generalization of the classical telegraph process which is an important model in the description of random motions on the real line. In its classical version, the telegraph process is a continuous-time random walk where a particle (the walker) moves with constant velocity which is reversed (or changed) at the arrival instants in a Poisson process (see \cite{Goldstein1951,Kac1974,Orsingher1990} and consult \cite{BogachevRatanov2011} for occupation time distributions of the telegraph process). 
A relativistic approach to remove drift terms is presented in \cite{BeghinNiedduOrsingher2001}.
In its classical (continuous-time) version the telegraph process (also called Cattaneo process) is Markovian reflecting the memoryless feature of the standard Poisson process.
Meanwhile, a large number of generalizations were introduced in the literature. Amongst them we mention here a model with occurrence of random velocities
\cite{StadjeZacks2004}, velocity reversals governed by a renewal process with IID Erlang-distributed interarrival times \cite{DiCrescendo2001}, and
further generalizations \cite{CrescenzoMartinucci2013} including semi-Markovian continuous time-fractional models and anomalous transport \cite{Masoliver2016,CompteMetzler1999,Gorskaetal2020} (consult also the references therein). 
To the best of our knowledge non-Markovian discrete-time versions of the telegraph process were so far not considered in the literature. The present paper is devoted to this subject.
\\[1mm]
The organization of our paper is as follows.
In Section \ref{fast-navigation} we introduce the SRW on the integer line with directed unit steps (simple walk) and discuss some relevant averages and their generating functions.
Large time asymptotic features are considered in Section \ref{large_times}.
In Section \ref{pertienent_GFs}
general formulae are derived for the sojourn probabilities on sites (the `spatial probability density function' or `propagator'). 
Section \ref{Bernoulli} is devoted to the `Bernoulli SRW', exhibiting geometrically distributed waiting times between the step reversals and standing out by the Markov property.
We derive explicit expressions for the expected position, mean square displacement and variance. 
In the diffusive limit (Section \ref{Bernoulli_telegraph}) the Bernoulli SRW converges to the classical telegraph process where the spatial probability density function solves a telegrapher's equation with drift. 
In Section \ref{continuous-time_limits} we consider a SRW with fat-tailed distributed ‘fractional Bernoulli’ waiting
times between the step reversals. The diffusive limit yields a fractional generalization of the telegraph process.
We analyze anomalous diffusive features
of the SRW (Section \ref{anomalous})
and take into account the `aging effect' as a sign of non-Markovianity where the discrete-time version of the so called `Aging Renewal Process' comes into play.
The Aging Renewal Process was introduced in the literature and applied to random walks for continuous times \cite{GodrecheLuck2001,Barkai2003,BarkaiCheng2003,SchulzBarkaiMetzler2014}. 
For the memoryless cases such as for Bernoulli SRW no aging effect occurs. In contrast, the `Sibuya SRW' with long waiting times between the step reversals exhibits a strong aging effect. For the Sibuya SRW we derive exact formulae for the variance in terms of discrete-time Prabhakar kernels exhibiting ballistic superdiffusive behavior in the large-time limit. The ballistic superdiffusive large time asymptotics is a general feature when the waiting times follow a fat-tailed PDF. In Section \ref{broad_densities} we explore diffusive features for two classes of SRWs where the waiting time densities have finite means.
The whole demonstration is accompanied by appendices where we deduce exact formulae for pertinent GFs for distributions of the `{\it Discrete-Time Aging Renewal Process}' (DTARP) (Appendices \ref{aging_discrete-time_renewal}, \ref{pertinent_DTARP}).
\\[1mm]
Finally, in Section \ref{examples} we introduce the 
`{\it continuous-time squirrel random walk}' (CTSRW) and consider the SRW subordinated to
a continuous-time renewal process
such as the 
time-fractional Poisson process. 
The CTSRW defines a
class of continuous-time random walks with different regimes of anomalous 
diffusion which we discuss by means of some examples.

\section{Squirrel random walk}
\subsection{General definition and preliminaries}
\label{fast-navigation}
We define the SRW as a discrete-time walk
characterized by the random variables $X_{t} \in \mathbb{Z}$, ($t \in \mathbb{N}_0$ denotes the integer-time)
with unit steps $\sigma_t=\{-1,1\}$ to the right or left direction as follows
\beq
\label{random_var}
X_t =  \sum_{r=1}^t \sigma_r , \hspace{1cm} X_0=0 , \hspace{1cm} t=1,2,\ldots
\eeq
where at $t=0$ the squirrel is sitting in the origin. 
A precise definition follows hereafter. We consider a discrete-time 
counting (renewal) process ${\cal N}(t) \in \mathbb{N}_0$ such that 
\cite{PachonPolitoRicciuti2021,MichelitschPolitoRiascos2021,TMM_FP_APR_Prab_Warsaw2021,ADTRW2021}
\beq
\label{discrete-time}
{\cal N}(t) = max(n\geq 0: J_n \leq t)  , \hspace{1cm} {\cal N}(0)=0 , \hspace{1cm} t =0,1,2\ldots
\eeq
with the arrival times (time instants of events, renewals) characterized by the random variables
\beq
\label{renewal_chain}
J_n=\sum_{j=1}^n \Delta t_j , \hspace{1cm} J_0=0.
\eeq
This is a discrete version of a strictly increasing
subordinator (see \cite{PachonPolitoRicciuti2021} for details and the references therein) with IID (independent and identically distributed) strictly positive integer increments $\Delta t_j  \in \mathbb{N} = \{1,2,\ldots\}$ (IID interarrival 
times or synonymously `waiting times' in the renewal interpretation). The increments are drawn from the discrete-time probability density function (PDF) of the form\footnote{We use the synonymous notations $\psi_t=\psi(t)$ and employ the term `probability density function' (PDF) or simply `waiting time density' for both discrete and continuous time cases.}
\beq
\label{first_arrival}
\mathbb{P}(\Delta t=k)= \psi_k , \hspace{1cm}  k=1,2,\ldots
\eeq
supported on positive integers $k\in \mathbb{N}$.
A brief remark on notations employed in this paper: In some sections (especially when additional indices come into play) we switch from index to function notations, for instance $f_r$ is then replaced by $f(r)$ for functions of $r$. To avoid any inconsistencies we clearly point out when we locally change notations.
\\[1mm]
The strictly increasing random process (\ref{renewal_chain}) is called a renewal chain and its inverse the discrete-time  (counting) renewal process (\ref{discrete-time}) represents the number of events (renewals) up to time $t$.
It is useful to introduce the generating function (GF) of the waiting time PDF,
\beq
\label{gen_fu}
\left\langle u^{\Delta t} \right\rangle= {\bar \psi}(u) = \sum_{t=1}^{\infty}\psi_tu^t , \hspace{1cm} |u| \leq 1,
\eeq
which fulfills ${\bar \psi}(u)\big|_{u=1}=1$ indicating normalization of (\ref{first_arrival}).
We introduce a random variable $Z_t \in \{0,1\}$ with $Z_t=1$ 
if there is an event (renewal) at instant $t$ and $Z_t=0$ at uneventful time instants. In the picture where the discrete-time renewal process is interpreted as a trial process, $Z_t=1$ indicates
a `success' (event) and $Z_t=0$ a `fail' (no event) in the trial performed at time $t$. We define the initial condition $Z_0=0$ (no `success' or event at $t=0$).
The counting variable (\ref{discrete-time}) can then be represented as
\beq
\label{repres}
{\cal N}(t) =\sum_{k=1}^t Z_k , \hspace{1cm} {\cal N}(0)=0.
\eeq
Now we connect the directed random walk (\ref{random_var}) with this counting process in the following way:
\begin{enumerate}
\item [(i)] At uneventful time instants $t$, i.e. $Z_t=0$,
the squirrel performs a unit step $\sigma_t=\sigma_{t-1}$ in the same direction as at $t-1$ where this holds for $t\geq 2$.
\item [(ii)] At arrival times $t$, i.e. $Z_t=1$, 
the step direction changes with respect to the previous step $\sigma_t=-\sigma_{t-1}$. 
\item [(iii)] 
We define that the first step is performed at $t=1$ in the direction
$\sigma_1=(1-2Z_1){\tilde \sigma}_0$ where ${\tilde \sigma}_0 \in \{-1,1\}$, i.e. the first step is in ${\tilde \sigma}_0$-direction if $t=1$ is uneventful
and in $-{\tilde \sigma}_0$-direction if $t=1$ is arrival time. 
The direction ${\tilde \sigma}_0$ can be thought as
either prescribed or randomly chosen.
\end{enumerate}
A little variant of the SRW is obtained when at each renewal time $J_n$ 
the sign reversal of the step is performed with a certain probability $p$ and with complementary probability $1-p$ the step direction remains unchanged. 
In other words the step direction is reversed according to (i)--(iii) at
arrival times of a new (composed) counting process ${\cal N}_B[{\cal N}(t)$] (${\cal N}_B(t)$ being a Bernoulli counting process \cite{PachonPolitoRicciuti2021} independent of ${\cal N}(t)$). Therefore, this variant refers also to the class of SRWs and does not define a further type of walk. 
\\[2mm]
With these considerations we can establish a simple recursion for the steps (see Eq. (\ref{random_var}))
 \beq
 \label{recursion_memroy}
 \sigma_{t} =(-1)^{Z_t}\sigma_{t-1} =(1-2Z_t)\sigma_{t-1}, \hspace{1cm} t\geq 2
 \eeq
with $\sigma_1=(-1)^{Z_1}{\tilde \sigma}_0$ and initial condition $\sigma_t\big|_{t=0}=\sigma_0=0$ (not equal to ${\tilde \sigma}_0$), i.e.\ no step at $t=0$ to fulfill initial condition $X_0=0$.
Hence, the increment at time $t$ can be represented by
\beq
\label{recrseive_val}
\sigma_t = {\tilde \sigma}_0[(-1)^{{\cal N}(t)} -\delta_{t0}] , \qquad t\geq 0,
\eeq
where $\delta_{rs}$ indicates the Kronecker symbol\footnote{We use the synonymous notation $\delta_{i,j}=\delta_{ij}$ for the Kronecker symbol.}.
The random variable (\ref{random_var}) then becomes
\beq
\label{ran_var}
X_t = X_{t-1} + {\tilde \sigma}_0(-1)^{{\cal N}(t)} ,\hspace{1cm} t =1,2,\ldots
\eeq
Now let us introduce the state probabilities 
$\mathbb{P}[{\cal N}(t)=n]=\Phi^{(n)}(t)$ denoting the probabilities 
for $n=0,1,2,\ldots$ arrivals within the discrete time interval $[0,t]$. We note that $\Phi^{(n)}(t)= 0$ for $n>t$ as a consequence that ${\cal N}(t) \leq t$ almost surely and of the initial condition
$\Phi^{(n)}(t)\big|_{t=0}=\delta_{n0}$. 
In order to compute the average position of the squirrel it is useful to consider the
generating function
\beq
\label{expect_increment}
\left\langle v^{{\cal N}(t)}\right\rangle  = {\cal P}(v,t) = 
\sum_{n=0}^t \mathbb{P}[{\cal N}(t)=n] v^n \qquad |v| \le 1.
\eeq
For a discrete-time counting process this GF is a
polynomial of degree $t$.
We called this polynomial in a former work `{\it state polynomial}' ${\cal P}(v,t)$ of the counting process \cite{ADTRW2021}. 
Considering ${\tilde \sigma}_0$ given, the mean increment at time $t$ writes
\beq
\label{expect_gen}
\begin{array}{clr}
\ds \langle \sigma_t\rangle &= \ds {\tilde \sigma}_0\left\langle (-1)^{{\cal N}(t)} -\delta_{t0}\right\rangle = 
{\tilde \sigma}_0 \sum_{n=0}^t \mathbb{P}({\cal N}(t)=n) [(-1)^n -\delta_{t0}]  &   \\ \\
 & =  \ds {\tilde \sigma_0}[{\cal P}(-1,t) - \delta_{t0}] &
 \end{array} \hspace{0.5cm} \ds (t\geq 0)
\eeq
with $\langle \sigma_t\rangle\big|_{t=0}=0$. Separating the even and odd
event numbers we get the probabilities that at time $t$ the step direction
is ${\tilde \sigma}_0$ i.e.\ for state $|+\rangle$) and $-{\tilde \sigma}_0$ for state $|-\rangle$,
where we come back later to this interpretation. Thus we get
\beq
\label{proba}
\begin{array}{clr}
\ds {\mathbb P}(\sigma_t{\tilde \sigma}_0=\sigma,t) &= \ds \frac{1}{2}\left\langle\left[1+(-1)^{{\cal N}(t)}\sigma\right]\right\rangle  &  \\ \\
 & = \ds \frac{1}{2}\left[1+\langle \sigma_t \rangle \sigma{\tilde \sigma}_0\right] &
 \end{array}  \hspace{1cm} \sigma = \pm 1, \hspace{1cm} t = 1,2,\ldots 
\eeq
which picks up in (\ref{expect_gen}) the terms with even $n$ for $\sigma=1$ and the terms with odd $n$ for $\sigma=-1$. Therefore,
\beq
\label{sigma_t}
\langle \sigma_t \rangle = {\tilde \sigma}_0 \sum_{\sigma =\pm 1} \sigma \, {\mathbb P}(\sigma_t{\tilde \sigma}_0=\sigma,t) , \hspace{1cm} t=1,2,\ldots
\eeq
It is then useful to evaluate the GF ${\bar \sigma}(u)$ of expression (\ref{expect_gen}) yielding
\beq
\label{mean_step_gen}
{\bar \sigma}(u) = \sum_{t=0}^{\infty} \langle \sigma_t\rangle u^t = 
{\tilde \sigma_0}\left[{\bar {\cal P}}(v,u)\bigg|_{v=-1}-1\right]
\eeq
where comes into play the GF of the state polynomial
\beq
\label{state_gen}
{\bar {\cal P}}(v,u) 
= \sum_{n=0}^{\infty}v^n{\bar \Phi}^{(n)}(u)  =\frac{1-{\bar \psi}(u)}{(1-u)[1-v{\bar \psi}(u)]} ,\hspace{1cm} |u|< 1, \hspace{1cm} |v| \leq 1
\eeq
in which we used the GF of the state probabilities
\cite{ADTRW2021,MichelitschPolitoRiascos2021}
\beq
\label{state_genfu}
{\bar \Phi}^{(n)}(u) =\sum_{t=0}^{\infty}\Phi^{(n)}(t)u^t=\frac{1-{\bar \psi}(u)}{1-u} ({\bar \psi}(u))^n.
\eeq
One gets, for the expected position, 
\beq
\label{mean_displaement}
\langle X_t \rangle = \sum_{r=0}^t\langle \sigma_r \rangle =
{\tilde \sigma}_0\left(\sum_{r=0}^t \left[{\cal P}(-1,r) -\delta_{r0}\right]\right)
\eeq
which has the GF
\beq
\label{gen_fu_of_position}
{\bar X}^{(1)}(u) = \sum_{t=0}^{\infty} u^t \langle X_t \rangle =
\frac{{\bar \sigma}(u)}{1-u} =
 \frac{[1-{\bar \psi}(u)]{\tilde \sigma}_0 }{(1-u)^2[1+{\bar \psi}(u)]}-\frac{{\tilde \sigma}_0}{1-u}
\eeq
with ${\bar X}^{(1)}(u)\big|_{u=0} = \langle X_0\rangle = 0$.
\\[2mm]
Consider now a sample path of the SRW up to time $t$ with the renewal chain (\ref{renewal_chain}) $J_1=\Delta t_1$,\, $J_2=\Delta t_1+\Delta t_2$, \,\ldots,\, $J_{{\cal N}(t)}=\Delta t_1+\ldots +\Delta t_{{\cal N}(t)}$. Per construction no step is performed at $t=0$ followed by
$\Delta t_1-1$ steps in ${\tilde \sigma}_0$-direction, $\Delta t_2$ steps in $-{\tilde \sigma}_0$-direction at the instants $\{J_1,\ldots, J_2-1\}$, $\Delta t_{{\cal N}(t)}$ steps in $(-1)^{{\cal  N}(t)-1}{\tilde \sigma}_0$-direction at the instants
$\{J_{{\cal N}(t)-1},\ldots, J_{{\cal N}(t)}-1\}$, and finally with $ t-J_{{\cal N}(t)}+1 \geq 1$ steps in
direction $(-1)^{{\cal N}(t)}{\tilde \sigma}_0$ at instants
$\{J_{{\cal N}(t)},\ldots, t\}$. This consideration leads to the representation
\beq
\label{repres1}
\begin{array}{clr} 
X_t  & =  \ds  {\tilde \sigma}_0\left(-1+\Delta t_1-\Delta t_2+\Delta t_3- \ldots +(-1)^{{\cal N}(t)-1}\Delta t_{{\cal N}(t)} + (-1)^{{\cal N}(t)}[t-J_{{\cal N}(t)}+1]\right)  & \\ \\
 &  = \ds  X_t^{(+)}-X_t^{(-)} & 
 \end{array}
\eeq
with $X_0=0$
where $|X_t| \leq t$ and time $t={\tilde \sigma}_0[X_t^{(+)}+X_t^{(-)}]$ counts the total number of steps made up to $t$. In this relation ${\tilde \sigma}_0X_t^{(+)}= \Delta t_1-1+\Delta t_3+\ldots$ indicates the number of steps in ${\tilde \sigma}_0$-direction and
${\tilde \sigma}_0 X_t^{(-)}= \Delta t_2+\Delta t_4+\ldots $ in the opposite direction up to time $t$. 
An interesting interpretation of the SRW is therefore that of
a two-state system (or connected graph of two nodes) where the squirrel is in state $|+\rangle$ (i.e. ${\cal N}(t)$ is even) at time instants when the step is in ${\tilde \sigma}_0$-direction and in state $\langle -|$ (${\cal N}(t)$ is odd) when a step is made in the opposite direction. The random variables
${\tilde \sigma}_0X_t^{(\pm)}$ can then be conceived as occupation times (sojourn times), i.e. the number of time units the squirrel has spent in these states during the interval $[0,t]$.
\\[1mm]
For the following we introduce the propagator $P(x,t)=: \mathbb{P}(X_t=x)$ which indicates the probability that the squirrel at time $t \in \mathbb{N}_0$ 
is present on the site $x\in \mathbb{Z}$. The propagator has the Fourier representation
\beq
\label{occupation_prob}
P(x,t)= \langle \, \delta_{x,X_t} \, \rangle =\left\langle \, \frac{1}{2\pi}\int_{-\pi}^{\pi} e^{i\kappa (x-X_t)}{\rm d}\kappa  \, \right\rangle = \frac{1}{2\pi}\int_{-\pi}^{\pi}e^{i\kappa x}\langle \, e^{-i\kappa X_t} \, \rangle {\rm d}\kappa , \hspace{1cm} (x \in \mathbb{Z}, \hspace{0.5cm} t \in \mathbb{N}_0)
\eeq
and fulfills initial condition $P(x,t)\big|_{t=0} = \delta_{x,0}$ as a consequence of $X_0=0$ and with the characteristic function
\beq
\label{Fourier_SRW_GF}
{\cal P}_{\kappa}(t) = \langle \, e^{-i\kappa X_t} \, \rangle .
\eeq
In section \ref{pertienent_GFs} we derive the GF of (\ref{Fourier_SRW_GF}) in explicit form (see Eqs. (\ref{sum_up_yields_occupation_GF}), (\ref{charfu_SRW})).
\subsection{Large time asymptotics of the expected position}
\label{large_times}
Before doing so, let us consider the large-time asymptotics of the expected position $\langle X_t \rangle$ which can be extracted
from (\ref{gen_fu_of_position}) for $u \to 1$ by virtue of Tauberian arguments. 
We notice the general asymptotic feature \cite{ADTRW2021}
\beq
\label{gen_feat}
{\bar \psi}(u) =1-A_{\mu}(1-u)^{\mu}+o[(1-u)^{\mu}] , \hspace{1cm} \mu \in (0,1], \hspace{0.5cm} (u \to 1)
\eeq
with $A_{\mu}>0$ and the admissible index range $\mu \in (0,1]$.
For $\mu \in (0,1)$ the waiting time density is fat-tailed (FT) 
and scales as $t^{-\mu-1}$ for large times,
where the waiting time has infinite mean $\frac{d}{du}{\bar \psi}(u)|_{u=1} \to \infty$ (see \cite{PachonPolitoRicciuti2021,MichelitschPolitoRiascos2021} for some details).
The class with $\mu=1$ contains two important subclasses: 
(i) Waiting time densities 
where
all moments are finite implying that
$\frac{d^r}{du^r}{\bar \psi}(u)\big|_{u=1} < \infty$ ($\forall r \in \mathbb{N}_0$) where $A_1=\frac{d}{du}{\bar \psi}(u)\big|_{u=1} =\sum_{t=1}^{\infty} \psi_t t \geq 1$
(light-tailed - `narrow' waiting time densities falling off at least geometrically). (ii)
Broad waiting time densities (not light-tailed) where the mean waiting time is finite but the variance infinite. This suggests the asymptotic expansion ($u\to 1-$)
\beq
\label{infvar}
{\bar \psi}(u) = 1-A_1(1-u)+B_{\lambda}(1-u)^{\lambda} + o[(1-u)^{\lambda}] 
\eeq
where $B_{\lambda} > 0$ and $1<\lambda \leq 2$. The subclass (i) is contained for $\lambda=2$.
The class (ii) (broad densities) has non-integer $\lambda \in (1,2)$. 
For the rest of the paper we classify waiting-time densities as follows. We call the subclass (i) with $\mu=1$, $\lambda=2$ `narrow' or `light-tailed' (LT) densities, the subclass (ii) which has $\mu=1$, $1<\lambda<2$ we call `broad' densities and the class with $0<\mu<1$ `fat-tailed' (FT) densities.
For continuous times having the same type of large time behavior as for discrete times, the subclass (ii) has been first considered in \cite{GodrecheLuck2001} and see also \cite{SingGorskaSandev2022} for an application to stochastic resetting. The large time behavior of a broad waiting time density is then governed by the power-law scaling 
(we use notation ``$\sim$'' for asymptotic equality)
\beq
\label{power-broad}
\psi(t) \sim B_{\lambda} \frac{\Gamma(t-\lambda)}{\Gamma(t+1)\Gamma(-\lambda)} \sim B_{\lambda} \frac{t^{-\lambda-1}}{\Gamma(-\lambda)} \to 0 , \hspace{1cm} \lambda \in (1,2)
\eeq
where $-3 < -\lambda-1 <-2$ having hence a lighter tail as a FT density.
Be reminded that $\Gamma(-\lambda)=\frac{\Gamma(2-\lambda)}{\lambda(\lambda-1)} >0$ as $\lambda \in (1,2)$. In the light-tailed limit $\lambda \to 2-$ (\ref{power-broad})
converges to $B_2\frac{d^2}{dt^2}\delta(t)=0$ thus the tail dies out reflecting the rapid decay of LT densities.
The GF (\ref{gen_fu_of_position}) of the expected position yields 
\beq
\label{genera_as_uto_one}
\ds {\bar X}^{(1)}(u) \sim  \left\{\begin{array}{ll} \ds {\tilde \sigma}_0 \frac{A_{\mu}}{2} (1-u)^{\mu-2} +o[(1-u)^{\mu-2}],  & \mu \in (0,1)  \\ \\
 \ds {\tilde \sigma}_0 \left(\frac{A_{1}}{2} -1\right) (1-u)^{-1} - {\tilde \sigma}_0\frac{B_{\lambda}}{2}(1-u)^{\lambda-2} +o[(1-u)^{\lambda-2}], & \mu=1 , \hspace{0.5cm}
 \lambda \in (1,2]\end{array}\right. 
\eeq
which yields for $\mu \in (0,1)$ a power-law escape for large times:
\beq
\label{asymp_pos__large_t}
\ds \langle X_t \rangle  \sim \left\{\begin{array}{ll} \ds \frac{{\tilde \sigma}_0A_{\mu}}{2}\frac{\Gamma(t+2-\mu)}{\Gamma(2-\mu)\Gamma(t+1)} \sim \frac{{\tilde \sigma}_0A_{\mu}}{2}\frac{t^{1-\mu}}{\Gamma(2-\mu)} \to \infty , \hspace{1cm}  \mu \in (0,1) & \\ \\  
\ds \frac{{\tilde \sigma}_0}{2} (A_1-2) - {\tilde \sigma}_0\frac{B_{\lambda}}{2}\frac{\Gamma(t+2-\lambda)}{\Gamma(2-\lambda)\Gamma(t+1)} \sim 
       \frac{{\tilde \sigma}_0}{2} (A_1-2) - 
       {\tilde \sigma}_0\frac{B_{\lambda}}{2}\frac{t^{1-\lambda}}{\Gamma(2-\lambda)} \to  \frac{{\tilde \sigma}_0}{2} (A_1-2)
   ,  &  \mu=1. \end{array}\right.
\eeq
For $\mu=1$ the squirrel in the average remains trapped where
the expected position approaches the limit value $\frac{{\tilde \sigma}_0}{2} (A_1-2)$ either (i) at least geometrically or (ii) slower with a $t^{1-\lambda}$ power-law. 
\\[1mm]
As $|\langle X_t\rangle | \leq t$ is bounded, the escape of the squirrel cannot be faster than $t$. The only class where the squirrel escapes to infinity is the FT class $0<\mu<1$ where the escape is sublinear with a $t^{1-\mu}$-law and in the direction of ${\tilde \sigma}_0$. The power-law escape behavior can be explained by long waiting times between the step direction reversals and can be interpreted as a fractal scaling of the expected position with respect 
to time (number of steps). 
\\[1mm]
As said $\mu=1$ represents the class with finite mean waiting times
$A_1=\frac{d}{du}{\bar \psi}_1(u)\big|_{u=1} < \infty$ where $A_1 \geq 1$ as $\psi_t$ is supported on $\mathbb{N}$. The special case if $A_1=1$ which corresponds to the minimum waiting time $\Delta t=1$ a.s. (i.e. $\psi(t)=\delta_{t1}$), represented by the trivial counting process ${\cal N}(t)=t$, yields an oscillatory and deterministic motion $\langle X_t \rangle = X_t=\frac{{\tilde \sigma}_0}{2}[(-1)^t -1]$ with GF ${\bar X}^{(1)}_1(u)=\frac{-u{\tilde \sigma}_0}{(1+u)(1-u)} \sim -\frac{\tilde \sigma_0}{2(1-u)}$ ($u\to 1-$) with the large-time average 
$\langle X_t \rangle \sim -\frac{{\tilde \sigma}_0}{2}$ of (\ref{asymp_pos__large_t}). This oscillatory motion
also occurs for $p=1$ in the `Bernoulli SRW' considered in Section \ref{Bernoulli}.
For $A_1=2$ the walk is asymptotically unbiased with $ \langle X_t \rangle = 0$ for large $t$ where in the long-time average the direction of each second step is reversed behaving as the symmetric Bernoulli SRW with $p=1/2$.
\subsection{Limit of infinite waiting times $\mu \to 0+$: `frozen limit'}
\label{frozen_limit}
It is worthy to consider here the limit of infinite waiting times characterized by the limit $\mu=\epsilon \to 0+$ more closely. Without loss of generality we assume $A_{\mu}=1$ corresponding to Sibuya distributed waiting times (considered subsequently in details, see (\ref{Sibuya}), (\ref{Sibuya_PDF})). In this limit for finite $t$, the squirrel is trapped in the state $|+\rangle$ (corresponding to ${\cal N}(t)=0$) performing steps solely in ${\tilde \sigma_0}$-direction (a.s.).
In particular, we are interested in the time range $0 \leq t < t_{\epsilon} \to \infty$ (as $\epsilon \to 0+$) for which the survival probability (which scales for large $t$ as 
 $\mathbb{P}[{\cal N}_{\epsilon}(t)=0] \sim \frac{1}{\Gamma(1-\epsilon)} t^{-\epsilon} \sim 1-$ \cite{ADTRW2021}) remains close to one:
$(t_{\epsilon})^{\epsilon} \sim 1+$, i.e.
$\epsilon \ln(t_{\epsilon}) \to 0+$. By assuming the scaling
$\ln(t_{\epsilon}) \sim \epsilon^{-\delta_1} (\gg 1)$
with $\delta_1 \in (0,1)$, this defines for $\epsilon \to 0+$ an infinitely large time interval without events. Therefore,
\beq
\label{limitsmallmu}
\langle X_t \rangle = {\tilde \sigma}_0 t ,\hspace{1cm} 0\leq t \leq t_{\epsilon} \approx \exp{[\epsilon^{-\delta_1}]} \to \infty , \hspace{0.5cm}  \mu =  \epsilon \to 0+ 
\eeq
i.e. the deterministic strict walk without step reversals emerges for any finite time $t$. We point out that the strict walk (\ref{limitsmallmu}) does not coincide with (\ref{asymp_pos__large_t}) for $\mu=0+$ where the latter is reduced by a factor of $1/2$. We will consider this issue closely.
We can choose $\epsilon$ sufficiently small
that for any finite fixed time $t < \exp{(\epsilon^{-\delta_1})}$ the squirrel is trapped in state $|+\rangle$ making solely steps in ${\tilde \sigma}_0$-direction. 
The GF of (\ref{limitsmallmu}) is ${\sigma}_0u/(1-u)^2$ and corresponds to the (forbidden) value $\mu=0$ by setting the (forbidden) waiting-time GF ${\bar \psi}(u)=0$ in (\ref{gen_fu_of_position}).
\\[1mm]
On the other hand a further time scale exists where $t$ is chosen sufficiently large that 
$t^{\epsilon} \gg 1$ for any given $0<\epsilon \ll 1$, where the survival probability is considerably reduced and many reversals of step directions have occurred. Clearly this holds true for times $t > \exp{(\epsilon^{-\delta_2})} \to \infty$ with $\delta_2 > 1$ when $\epsilon \to 0+$. In this time range
the expected position then turns into the large time power-law (\ref{asymp_pos__large_t}), namely
\beq
\label{real_large_times}
\langle X_t \rangle \sim \frac{{\tilde \sigma}_0}{2}t^{1-\epsilon} , \hspace{1cm} t > \exp{[\epsilon^{-\delta_2}]} \gg \exp{[\epsilon^{-\delta_1}]} \to \infty , \hspace{0.5cm} \delta_2>1 , \hspace{0.5cm} \mu=\epsilon \to 0+
\eeq
where $\Gamma(2-\epsilon) \to 1$. The smaller $\mu$ the longer it takes for the squirrel to reach the large time power-law
(\ref{real_large_times}) and it takes infinitely long in the frozen limit $\mu \to 0+$.
\\[1mm]
With a similar consideration one can see for geometrically distributed waiting times (see \cite{ADTRW2021} for details) that the deterministic walk (frozen regime) (\ref{limitsmallmu}) emerges for $p \to 0+$ ($p$ Bernoulli probability of step reversal at instant $t$) 
for $0 \leq t < p^{-\delta_1} \to \infty$ ($\delta_1 \in (0,1)$) where the Bernoulli PDF 
$\psi_B(t) = p(1-p)^{t-1} \sim p 
(1-p)^{p^{-\delta_1}} \sim p \to 0+$
and the survival probability is close to one
$(1-p)^t > (1-p)^{p^{-\delta_1}}\sim \exp{[\,-p^{(1-\delta_1)}\,]} \to 1-$.
We consider the `Bernoulli SRW' thoroughly in Section \ref{Bernoulli}.
\section{Propagator and related generating functions}
\label{pertienent_GFs}
In this section we derive GFs defining the propagator of the SRW, i.e.\ the probabilities that the squirrel at time instant $t \in \mathbb{N}_0$ is present on site $x \in \mathbb{Z}$.
We use in the following and throughout the paper for expected values the notation
\beq
\label{notation_average}
\left\langle\, f(\Delta t) \, \right\rangle = \sum_{r=1}^{\infty}\psi(r)f(r)
\eeq
for the averaging over $\Delta t_j$ being IID copies of $\Delta t \in \mathbb{N}$ with PDF $\psi(t)$ defined in (\ref{first_arrival}). A simple example for (\ref{notation_average}) is the expected value
of the Kronecker symbol
$\langle \, \delta_{k,\Delta t}\, \rangle =\psi(k)$ recovering the waiting-time PDF.
For our convenience we introduce the following discrete step function defined on integers
\beq
\label{gen_Theta}
\Theta(a,r,b) = \Theta(r-a)-\Theta(r-b) = 
\left\{\begin{array}{ll} 1 &\hspace{1cm} \mathrm{for} \hspace{0.25cm} a \leq r \leq b-1
\\ 0 & \hspace{1cm} \mathrm{otherwise}
\end{array} \hspace{0.5cm} a < b,\hspace{0.5cm} (a,b, r \in \mathbb{N}_0)
\right.
\eeq
where we used the discrete Heaviside function
\beq
\label{Heaviside_step}
\Theta(k) = \left\{\begin{array}{clr} \ds  1 , & &\hspace{1cm}  \ds k\geq 0 \\ 
\ds 0, & & \ds \hspace{1cm} k < 0 \end{array}\right. \hspace{1cm} (k\in \mathbb{Z}).
\eeq
Note that $\Theta(0)=1$. The convenience of the step function (\ref{gen_Theta}) for our evaluations is that
$\Theta(J_n,t,J_{n+1})=1$ only for the $\Delta t_{n+1}$ time instants $J_n\leq t \leq J_{n+1}-1$, i.e.\ when ${\cal N}(t)=n$ and $\Theta(J_n,t,J_{n+1})=0$ for ${\cal N}(t)\neq n$.
Now, consider
\beq
\label{meanmultivariate}
f(t,\tau,\zeta_1,\ldots,\zeta_n) =  \langle \, \Theta(J_n,t,J_{n+1})\delta_{\tau,t-J_n} \zeta_1^{\Delta t_1}\ldots\zeta_n^{\Delta t_n} \, \rangle
, \hspace{1cm} (0< |\zeta_j|\leq 1, \hspace{0.5cm} \tau,\, t \in \mathbb{N}_0)
\eeq
and evaluate its double GF ($t\leftrightarrow u, \tau \leftrightarrow w$)
\begin{equation*}
{\bar f}(u,w,\zeta_1,\ldots,\zeta_n)= \sum_{k=0}^{\infty}\sum_{s=0}^{\infty}f(s,k,\zeta_1,\ldots,\zeta_n) w^{k}u^s,
\end{equation*}
which yields
\beq
\label{GF_multivariante}
\begin{array}{clr}
\ds {\bar f}(u,w,\zeta_1,\ldots,\zeta_n) & = \ds \langle \, \zeta_1^{\Delta t_1}\ldots\zeta_n^{\Delta t_n} w^{-J_n} \sum_{s=J_n}^{J_{n+1}-1} u^{s} w^{s} \, \rangle &\\ \\  & =
\ds  \langle\,(u\zeta_1)^{\Delta t_1}\ \rangle\ldots \langle\,  (u\zeta_1)^{\Delta t_n}\, \rangle  \frac{1- \langle\,(uw)^{\Delta t_{n+1}}\,\rangle}{1-uw}  \\ \\ & = \ds {\bar \psi}(u\zeta_1)\ldots {\bar \psi}(u\zeta_n) \, \frac{1-{\bar \psi}(uw)}{1-uw} &
\end{array}
\eeq
where we used the IID feature of the $\Delta t_j$ and for $n=0$ the empty product has to be read as equal to one.
Clearly, for $\zeta_{\ell}=1$ and $w=1$ this recovers GF (\ref{state_genfu}) of the state probabilities
$\mathbb{P}({\cal N}(t)=n) = \langle \, \Theta(J_n,t,J_{n+1}) \,\rangle$.
\\[2mm]
Now consider a generalization of the SRW  (\ref{repres1}) with directed
steps of prescribed sizes $a_n$ (instead of unit steps) for the $\Delta t_n$ time instants within 
$t \in [J_{n-1},J_n-1]$ ($n\geq 1$). The SRW we have considered so far is a special case with
$a_n= (-1)^{n-1}{\tilde \sigma}_0$.
A sample path of the generalized SRW is 
\beq
\label{gen_SRW}
(X_t)_{\{a_{\ell}\}} = a_1(\Delta t_1-1)+a_2\Delta t_2 +\ldots + a_{n} \Delta t_n + a_{n+1}(t-J_n+1) , \hspace{1cm} (n={\cal N}(t))
\eeq
with $(X_0)_{\{a_{\ell}\}}=0$.
The characteristic function of this walk $\langle \, e^{-i\kappa (X_t)_{\{a_{\ell}\}}}\, \rangle$ can be derived from the expected value 
\beq
\label{char_fuSRW}
\begin{array}{clr}
\ds g(t,\zeta_1,\zeta_2, \ldots,\zeta_n;\zeta_{n+1}) & = \ds  \langle \, 
\Theta(J_n,t,J_{n+1}) \zeta_1^{\Delta t_1-1}\ldots\zeta_n^{\Delta t_n}\zeta_{n+1}^{t-J_n+1}\,
\rangle & \\ \\  & = \ds  \zeta_1^{-1}\zeta_{n+1} \,\, \langle \, 
\Theta(J_n,t,J_{n+1}) \zeta_1^{\Delta t_1}\ldots\zeta_n^{\Delta t_n}\zeta_{n+1}^{t-J_n}\,
\rangle & 
\end{array}
\eeq
for $\zeta_j=e^{-i\kappa a_j}$. Its GF yields
\beq
\label{GF_n}
\begin{array}{clr}
\ds {\bar g}(u,\zeta_1,\ldots,\zeta_n;\zeta_{n+1}) & = \ds \zeta_1^{-1}\zeta_{n+1} {\bar f}(u,w,\zeta_1,\ldots,\zeta_n)\big|_{w=\zeta_{n+1}} & \\ \\
 & = \ds \zeta_1^{-1}\zeta_{n+1}  \, \frac{1-{\bar \psi}(u\zeta_{n+1})}{1-u\zeta_{n+1}} \, {\bar \psi}(u\zeta_1)\ldots {\bar \psi}(u\zeta_n) , \hspace{1.5cm} n =1,2,\ldots
 \end{array}
\eeq
and for $n=0$ we have ${\bar g}(u;\zeta_1) = \frac{1-{\bar \psi}(u\zeta_{1})}{1-u\zeta_{1}}$.
Summing up the ${\bar g}(\cdot;\zeta_{n+1})$ over $n$ yields
the GF of the characteristic function of the propagator for the generalized SRW. 
The SRW (\ref{repres1}) is represented by the case
$\zeta_{2\ell}= \zeta_2=e^{i\kappa{\tilde \sigma}_0}$ ($a_2=-{\tilde \sigma}_0$)
and $\zeta_{2\ell+1}=\zeta_1=e^{-i\kappa{\tilde\sigma}_0}$ ($a_1={\tilde \sigma}_0$). Then (\ref{GF_n}) yields
\beq
\label{special-SRW-case}
\ds {\bar g}_n(u,\zeta_1,\zeta_2)= \left\{\begin{array}{llr} \ds \frac{1-{\bar \psi}(u\zeta_1)}{1-u\zeta_1}\left[{\bar \psi}(u\zeta_1){\bar \psi}(u\zeta_2)\right]^{\ell}  , & n=2\ell & \\ \\ \ds \zeta_1^{-1}\zeta_2
\frac{1-{\bar \psi}(u\zeta_2)}{1-u\zeta_2}
{\bar \psi}(u\zeta_1)\left[{\bar \psi}(u\zeta_1){\bar \psi}(u\zeta_2)\right]^{\ell} , & n=2\ell+1. &
\end{array}\right.  (\ell=0,1,2,\ldots)
\eeq
Summing over $n$ takes us to the GF of the
characteristic function (\ref{Fourier_SRW_GF}):
\beq
\label{sum_up_yields_occupation_GF}
\begin{array}{clr}
\ds {\bar g}(u,\zeta_1,\zeta_2) & = \ds  \sum_{\ell=0}^{\infty}\left[{\bar g}_{2\ell}(u,\zeta_1,\zeta_2)+{\bar g}_{2\ell+1}(u,\zeta_1,\zeta_2)\right]  & \\ \\
& = \ds \frac{1}{1-{\bar \psi}(u\zeta_1){\bar \psi}(u\zeta_2)}\left(\frac{1-{\bar \psi}(u\zeta_1)}{1-u\zeta_1}+ \zeta_1^{-1}\zeta_2
\frac{1-{\bar \psi}(u\zeta_2)}{1-u\zeta_2}
{\bar \psi}(u\zeta_1)\right). & 
\end{array}
\eeq
Further,
\beq
\label{charfu_SRW}
{\bar P}_{\kappa}(u)= \left\langle\, \sum_{t=0}^{\infty} u^t e^{-i\kappa X_t} \,\right\rangle =
{\bar g}(u,e^{-i\kappa{\tilde \sigma}_0},e^{i\kappa{\tilde \sigma}_0}) , \hspace{1cm} \kappa \in (-\pi,\pi).
\eeq
For $u$ real ${\bar P}_{\kappa}(u) \in \mathbb{C}$, which means that the SRW is in general biased.
Let us now briefly check some necessary properties.
For $u=0$ we observe ${\bar g}(0,\zeta_1,\zeta_2)=1$
to fulfill in (\ref{occupation_prob}) the initial condition ${\bar P}(x,0)=\delta_{x,0}$.
Then, for $\zeta_1=\zeta_2=\zeta$ we have
$g(t,\zeta,\zeta) = \langle \zeta^t \rangle =\zeta^t$ in agreement with (\ref{sum_up_yields_occupation_GF}) to yield
\beq
\label{gxixi}
{\bar g}(u,\zeta,\zeta)= \frac{1}{1-\zeta u},
\eeq
thus confirming for $\zeta=1$ the normalization of the propagator.
\\[1mm]
It is worthy of mention that (\ref{GF_n})--(\ref{charfu_SRW}) have a rich field of applications. For instance, the propagators 
of the sojourn times (`occupation times') ${\tilde \sigma}_0X_t^{\pm}$ in the states $|\pm\rangle$ can be easily derived:
\beq
\label{occupation_times}
\begin{array}{clr}
\ds {\bar P}_{\kappa}^+(u) & = \ds \left\langle \, \sum_{t=0}^{\infty}u^t  e^{-i\kappa X^{+}_t} \,\right\rangle = {\bar g}(u,e^{-i\kappa {\tilde \sigma}_0},1) & \\ \\
\ds {\bar P}_{\kappa}^-(u)  & = \ds  \left\langle\, \sum_{t=0}^{\infty}u^t e^{-i\kappa X_t^{-}} \,\right\rangle = {\bar g}(u,1,e^{i\kappa {\tilde \sigma}_0}). &
\end{array}
\eeq
More complex cases of generalized SRW type (\ref{gen_SRW}) can be analyzed with this approach.
\\[2mm]
We prove now the following important feature:
\beq
\label{affirm}
\begin{array}{clr}
\ds \left(\frac{\partial}{\partial \zeta_1} - \frac{\partial}{\partial \zeta_2}\right)
H(t,\zeta_1,\zeta_2)\big|_{\zeta_1=\zeta_2=1} &  = \ds \sum_{r=0}^t{\cal P}(-1,r) = \langle\, \sum_{r=0}^t  (-1)^{{\cal N}(r)} \, \rangle  & \\ \\
& = \ds 1 + {\tilde \sigma}_0 \langle \, X^{+}_t-X^{-}_t\, \rangle & 
\end{array}
\eeq
(with the state polynomial ${\cal P}(v,r)$ defined in (\ref{expect_increment})) where we introduced
\beq
\label{axiliary_walk}
H(t,\zeta_1,\zeta_2) = \zeta_1 g(t,\zeta_1,\zeta_2) 
\eeq
corresponding to the auxiliary walk ${\tilde \sigma}_0X_t+1$ with a step in
positive direction at the (uneventful) initial time $t=0$ (a.s.). 
Now, with ${\bar H}(u,\zeta_1,\zeta_2) = \zeta_1 {\bar g}(u,\zeta_1,\zeta_2)$ 
we arrive at
\beq
\label{proof-of-above}
\begin{array}{clr}
\ds \left(\frac{\partial}{\partial \zeta_1} - \frac{\partial}{\partial \zeta_2}\right)
{\bar H}(u,\zeta_1,\zeta_2)\big|_{\zeta_1=\zeta_2=1} & = \ds \frac{1}{1+{\bar \psi}(u)}\left(\frac{u}{1-u}\frac{d {\bar \psi}(u)}{du} + \frac{d}{d\zeta}\left[\zeta\frac{1-{\bar \psi}(u\zeta)}{1-u\zeta}\right]\Big|_{\zeta=1}\right)   & \\ \\
 & = \ds \frac{1-{\bar \psi}(u)}{(1-u)^2[1+{\bar \psi}(u)]} = \frac{{\bar {\cal P}}(-1,u)}{1-u} & 
\end{array}
\eeq
which indeed is the GF of (\ref{affirm}), concluding the calculation.
It is only a little step to establish the connection with the GF of the expected position
\beq
\label{SRW_first_moment}
\begin{array}{clr}
\ds {\tilde \sigma}_0{\bar X}^{(1)}(u) & = \ds \left(\frac{\partial}{\partial \zeta_1} - \frac{\partial}{\partial \zeta_2}\right)\frac{{\bar H}(u,\zeta_1,\zeta_2)}{\zeta_1}\bigg|_{\zeta_1=\zeta_2=1} & \\ \\ & = \ds  \frac{{\bar {\cal P}}(-1,u)}{1-u}-{\bar H}(u,1,1) = \frac{{\bar {\cal P}}(-1,u)}{1-u}-\frac{1}{1-u} &
\end{array}
\eeq
in agreement with (\ref{mean_displaement}), (\ref{gen_fu_of_position}). 
\section{Bernoulli SRW}
\label{Bernoulli}
We consider now in more detail the SRW where step directions turn at arrival times of a Bernoulli counting process ${\cal N}_B(t)$ thus inheriting the Markov property.
We call this walk `Bernoulli SRW'. Let $p$ be the probability
of a success (event or arrival) in a Bernoulli trial and recall its geometrically distributed waiting time density:
$\psi_B(t,p)=pq^{t-1}$ ($t \in \mathbb{N}$, $ p\in (0,1]$, $q=1-p$). The value $p$ should be strictly positive otherwise the probability mass is concentrated at infinity.

The Bernoulli SRW has the following interpretation: at each integer time instant the squirrel turns the step direction with probability $p$ and maintains it with complementary probability $q=1-p$.
The GF 
(\ref{mean_step_gen}) of the expected steps has then the form
\beq
\label{Bernoulli_genfu}
{\bar \sigma}_B(u) = {\tilde \sigma}_0 \frac{(1-2p)u}{1-u(1-2p)}
\eeq
with the expected value of the increment 
\beq
\label{expect_pos_Bernoulli}
\langle \sigma_t\rangle_B = {\tilde \sigma}_0 \left[\langle (-1)^{{\cal N}_B(t)}\rangle -\delta_{t0}\right]    = {\tilde \sigma}_0 \left[ (1-2p)^t -\delta_{t0}\right] , \hspace{1cm} t =0,1,2,\ldots
\eeq
fulfilling the initial condition $\langle \sigma_0\rangle_B=0$
and where $\langle \sigma_1\rangle_B = {\tilde \sigma}_0 (q-p)$
for the first step.
%with
%$\langle (-1)^{{\cal N}_B(t)}\rangle = (1-2p)^t$ ($t=0,1,2\ldots$).
In the limit $p \to 0+$ no event occurs
%(for a very long time) where
and the walk becomes deterministic with $\langle \sigma_B(t)\rangle ={\tilde \sigma}_0$ 
where the steps do not change the
direction. 
The case with $p=1$ recovers the deterministic trivial counting process ${\cal N}_B(t)=t$ where the squirrel changes (a.s.) at each time increment its step direction
with oscillatory behavior $\langle \sigma_B(t)\rangle =(-1)^t {\tilde \sigma}_0$ ($t>0$). For $p=\frac{1}{2}$ the walk is unbiased with $\langle \sigma_B(t)\rangle =0 $ where the squirrel in the average remains on the departure site.
The Markov property of the Bernoulli SRW is reflected by the following feature:
\beq
\label{memoryless_feature_Bernoulli}
\langle \sigma_{t_1+t_2}\rangle_B =  \langle \sigma_{t_1}\rangle_B  \langle \sigma_{t_2}\rangle_B , \hspace{1cm} t_j=1,2,\ldots
\eeq
For later use we consider the GF of the expected position
\beq
\label{later_use}
{\bar X}^{(1)}_B(u) = \frac{{\bar \sigma}_B(u)}{1-u} =  {\tilde \sigma}_0 \frac{(1-2p)u}{(1-u)[1-u(1-2p)]},
\eeq
which yields straight-forwardly
\beq
\label{expected_pos_Bern}
\langle X_t\rangle_B = {\tilde \sigma}_0 \frac{1-2p}{2p}[1-(1-2p)^t] ,\hspace{1cm} t = 0,1,2,\ldots
\eeq
and is consistent with the general asymptotic relation for light-tailed waiting time densities (\ref{asymp_pos__large_t}) (identify with $A_1=\frac{1}{p}$, the expected waiting time).
For large $t$, the squirrel is localized approaching geometrically the value $\langle X_{t\to \infty}\rangle_B = {\bar \sigma}_B(u)\big|_{u=1} = {\tilde \sigma}_0 \frac{1-2p}{2p}$
which is located for $p<\frac{1}{2}$ on the same side as ${\tilde \sigma}_0$ and for $p>\frac{1}{2}$ on the opposite side.  
For $p=\frac{1}{2}$ the SRW is unbiased and, in the average, the squirrel turns the direction at any second time instant and hence its expected position remains localized on the departure site $X_0=0$. It is straight-forward to see that the only counting process which generates $\langle X_t \rangle = \langle \sigma_t\rangle = 0$ for all $t$ is the symmetric Bernoulli process. Putting the GF (\ref{gen_fu_of_position}) equal to zero defines a condition for ${\bar \psi}(u)$ which yields the symmetric Bernoulli generating function
${\bar \psi}_B(u)=u/(2[1-\frac{u}{2}])$.
\\[2mm]
As said, if $p=1$ ($\psi_B(t)=\delta_{t1}$) the trivial (deterministic) counting process ${\cal N}(t)=t$ is recovered with the deterministic oscillatory motion $X_t=\frac{{\tilde \sigma}_0}{2}[(-1)^t -1]$.
\\[2mm]
On the other hand, in the limit $p\to 0+$ (no arrival ${\cal N}(t)=0$ a.s. corresponding to the frozen limit considered in Section \ref{frozen_limit}) the squirrel 
maintains the direction of the initial step ${\tilde \sigma_0}$.
This limit is contained in (\ref{expected_pos_Bern}) by applying de L'H\^opital's rule
\beq\label{delopital}
\langle X_t\rangle \sim {\tilde \sigma}_0\frac{d}{d\xi}[1-(1-\xi)^t]\big|_{\xi=0} = {\tilde \sigma}_0 t .
\eeq
This deterministic limit (where the waiting time between step reversals becomes infinite) represents the fastest possible escape from the departure site with $|\langle X_t\rangle |=t$.
\\[1mm]
We can directly evaluate the mean square displacement
\beq
\label{variance_bernoulli}
\begin{array}{clr}
 \langle X_t^2\rangle_B  & =\ds \sum_{r=1}^t \sum_{k=1}^t \langle \sigma_k\sigma_r\rangle =-t+2K_B(t)  = \ds -t + 2\sum_{r=1}^{t} 
\sum_{k=r}^{t} \langle (-1)^{{\cal N}_B(k-r)}\rangle & \\ \\ & =  \ds  -t +  \frac{2 {\tilde \sigma}_0}{1-2p}\sum_{r=1}^{t} \langle X_r \rangle_B & 
\end{array}  \ds (t=1,2,\ldots)
\eeq
and see Eqs. (\ref{Bernoulli_K}), (\ref{Ber_msdgenfu}) for the GF of $K_B(t)$.
To obtain this result directly it is convenient to use the Markovian property of Bernoulli, i.e. for $k\geq r$ the quantity
${\cal N}_B(t)-{\cal N}_B(r) = {\cal N}_B(t-r) $ ($t\geq r$) is itself an independent Bernoulli counting variable without the effect of `aging' which we consider extensively in the Appendices. 
Now with $\langle (-1)^{{\cal N}_B(k)}\rangle = (1-2p)^k$ (see Eq. (\ref{expect_pos_Bernoulli})) we arrive at 
\beq
\label{meanquare_bern}
\begin{array}{clr}
\ds \langle X_t^2\rangle_B & = \ds  -t+\frac{1}{p}\sum_{r=1}^{t}[1-(1-2p)^r]= \frac{1-p}{p}t-\frac{1-2p}{2p^2}\left[1-(1-2p)^t\right] & \\ \\
 & =\ds  \frac{(1-p)}{p}\, t -\frac{{\tilde \sigma}_0}{p}\langle X_t\rangle_B  & 
 \end{array}  t=1,2,\ldots
\eeq
with $\langle X_0^2\rangle_B=0$ reflecting the initial condition. For $t=1$ (\ref{meanquare_bern}) necessarily yields $\langle X_1^2\rangle_B = \sigma_1^2=1$. Then we get
for the variance
\beq
\label{variance}
{\cal V}_B(t) =  \langle X_t^2\rangle_B - (\langle X_t\rangle_B)^2 =  \frac{(1-p)}{p}\, t -\frac{{\tilde \sigma}_0}{p}\langle X_t\rangle_B - (\langle X_t\rangle_B)^2 , \hspace{1cm} p\neq 0
\eeq
where for large times we have linear (normal diffusive) behavior
${\cal V}_B(t) \sim \langle X_t^2\rangle_B \sim \frac{(1-p)}{p}\, t$.
For $p<\frac{1}{2}$ the linear increase is faster $\frac{(1-p)}{p}>1$ than for $p>\frac{1}{2}$ with $\frac{(1-p)}{p}< 1$.
\\[2mm]
For $p=1$ we have an oscillatory (non-fluctuating) deterministic motion $X_t=\frac{{\tilde \sigma}_0}{2}((-1)^t-1)$ with ${\cal V}(t) = 0$ without linear increase of the variance. The mean square position remains bounded by the oscillating behavior $\ds \langle X_t^2\rangle_B = \frac{1}{2}(1-(-1)^t)$ (being null for even $t$ including $t=0$ and $+1$ for $t$ odd).
\begin{figure*}[t]
\centerline{\includegraphics[width=0.6\textwidth]{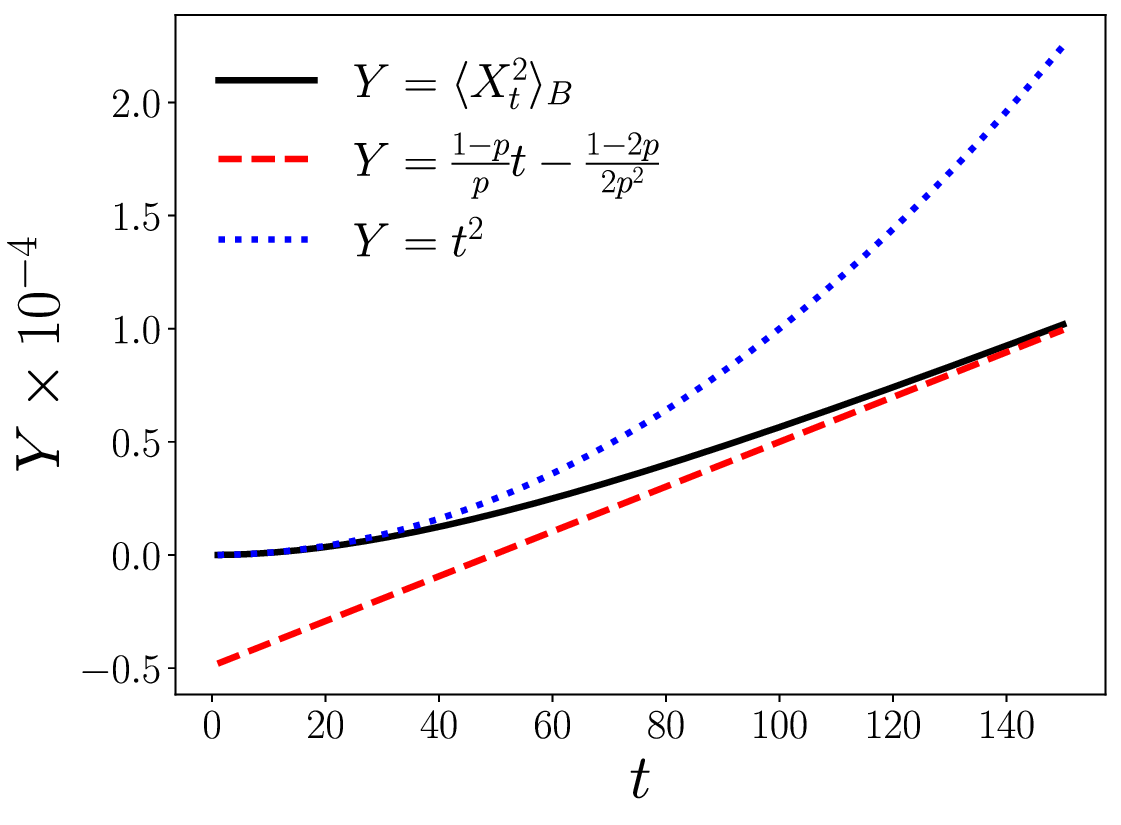}}
\vspace{-5mm}
\caption{\label{X2ber} Different curves $Y$ as a function of $t$. Mean square displacement $Y=\langle X_t^2\rangle_B$  (solid black curve) from Eq. (\ref{meanquare_bern}) in rescaled units for Bernoulli probability $p=0.01$.}
\end{figure*}
In the other deterministic limit $p\to 0+$ 
with $\langle X_t\rangle ={\tilde \sigma}_0 t $ we get by expanding (\ref{meanquare_bern}) with respect to $p$
(where some divergent terms cancel each other)
\beq
\label{expand_p}
\langle X_t^2 \rangle_B = t^2 , \hspace{1cm} (p \to 0+).
\eeq
Therefore, we get in this limit ${\cal V}(x)= t^2-({\tilde \sigma}_0 t)^2 \to 0 $ corresponding to the deterministic (non-fluctuating) motion 
without step reversals. 
We depict the mean square displacement for a small Bernoulli probability $p=0.01$ (rare change of step directions) in Fig. \ref{X2ber}. For small $t$ the behavior is nearby quadratic and close to the deterministic limit $p=0+$ (upper dotted curve). For larger times the linear asymptotics (lower dotted line) is geometrically approached.
\\[1mm]
Finally for $p=\frac{1}{2}$ the walk is unbiased with $\langle X_t\rangle=0$. Thus we have ${\cal V}_B(x)=\langle X_t^2\rangle_B = t $, corresponding to symmetric normal diffusion.
\section{Continuum limits of SRW}  
\label{Diffusive_telepgraph}
\subsection{Rescaled SRW}
In this part we consider the combined continuous space-time limit by simultaneously rescaling time and space units. In the rescaled SRW the directed steps $\Delta x=v_0h$
($v_0=|v_0|{\tilde \sigma}_0$ indicates the directed velocity independent of the time increment $h$) are performed at time instants $t\in\{0,h,2h,\ldots \}\in h
{\mathbb N}_0$. A sample path of the rescaled SRW is represented by (see (\ref{repres1}))
\beq
\label{rescaled_SRW}
\begin{array}{clr}
(X_t)_{hv_0} =v_0\left[-h+\Delta t_1-\Delta t_2+ \ldots +(-1)^{{\cal N}_h(t)-1}\Delta t_{{\cal N}_h(t)} + (-1)^{{\cal N}_h(t)}(t-J_{{\cal N}_h(t)}+h)\right]  & &
\end{array}
\eeq
with IID interarrival intervals $\Delta t_j=\{h,2h,\ldots\} \in h\mathbb{N}$ of the rescaled counting process
${\cal N}_h(t)=\max[n \in \mathbb{N}_0: J_n(h) \leq t \in h\mathbb{N}_0]$ where $J_n(h) \in h\mathbb{N}_0$ represents the rescaled renewal chain.
In order to obtain an existing diffusive limit it is necessary to rescale the time scale parameter in the waiting time density (consult \cite{PachonPolitoRicciuti2021,MichelitschPolitoRiascos2021,TMM_FP_APR_fractiona_fract} for extensive outlines and applications of `well-scaled' limits).
Whether the squirrel occupies a certain position on the lattice $x \in h|v_0| \mathbb{Z}$ can be expressed by the Kronecker symbol
$\delta_{\frac{x}{h|v_0|},\frac{(X_t)_{hv_0)}}{|v_0|h}}$.  
Therefore, considering
(\ref{occupation_prob}), the propagator writes 
\beq
\label{diffusive_limit-scaled}
\begin{array}{clr} 
\ds P_h(x,t,v_0) & = \ds \frac{1}{|v_0|h} \langle \, \delta_{\frac{x}{h|v_0|},\frac{(X_t)_{hv_0}}{|v_0|h}} \rangle & \\ \\ &  = \ds \left\langle\,  \frac{1}{2\pi h|v_0|}\int_{-\pi}^{\pi} e^{i\kappa \frac{(x-(X_t)_{hv_0})}{h|v_0|}}{\rm d}\kappa \,\right\rangle & \\ \\ & = \ds \left\langle\, \frac{1}{2\pi}\int_{-\frac{\pi}{h|v_0|}}^{\frac{\pi}{h|v_0|}} 
e^{ik(x-(X_t)_{hv_0})} {\rm d}k \,\right\rangle  & 
\end{array}
\eeq
where $P_h(x,t,v_0) h|v_0|$ denotes the probability to find the squirrel at time $t$ on site $x$.
The multiplier $1/{(|v_0|h)}$ comes into play as $P_h(x,t,v_0)$ is a spatial density (having units {$[length^{-1}]$}) attributed to interval  $[x,x+|v_0|h)$.
The rescaled propagator (\ref{diffusive_limit-scaled}) is normalized as
\beq
\label{normprop}
\sum_{r=-\infty}^{\infty}  P_h(rh|v_0| ,t,v_0) h|v_0| = \left\langle \,\sum_{r=-\infty}^{\infty}  \delta_{\ds r,(X_t)_{hv_0}/(|v_0|h)} \, \right\rangle    = 1.
\eeq
For $h\to {\rm d}t$ ($h|v_0| \to {\rm d}x$) the scaled SRW converges to the continuum limit: $(X_t)_{hv_0} \to X_{v_0,t}^{(c)} \in \mathbb{R}$ ($t\in \mathbb{R}^{+}$) with propagator
$P_h(x,t,v_0) \to P_c(x,t,v_0) = 
\langle \delta(x-X_{v_0,t}^{(c)})\rangle $ ($t \in \mathbb{R}^{+}$ and $x, X_{v_0,t}^{(c)} \in \mathbb{R}$) where $\delta(x-x')$ indicates infinite space Dirac's $\delta$-distribution.
In what follows we extensively use the features of Laplace transforms of causal functions and distributions (see Appendix \ref{Lapl_Trafo} and \cite{MichelitschRiascos2020b} for some details).
\subsection{Continuum limit of Bernoulli SRW to the telegraph process}
\label{Bernoulli_telegraph}
We now explore a `well-scaled' continuum limit for the Bernoulli SRW.
To this end we rescale the Bernoulli probability as $p=\xi_0h$ with the new time-scale constant $\xi_0>0$ of units $\sec^{-1}$ and independent of the time increment $h$.
The waiting time PDF of rescaled Bernoulli then is
$h^{-1}\psi_{B}(t/h,\xi_0h)= \xi_0(1-\xi_0h)^{t/h-1} \to \xi_0 e^{-\xi_0t}$ ($t \in h\mathbb{N}_0 \to \mathbb{R}^+$) and converging to the continuous-time exponential waiting-time density of the Poisson process. This reflects the fact that the Bernoulli process is a discrete version of the Poisson process both standing out by the Markov property.
Consider now the scaling limit of (\ref{sum_up_yields_occupation_GF}), (\ref{charfu_SRW}) for the scaled Bernoulli SRW with
$\zeta_1=e^{-ihkv_0}$, $\zeta_2=e^{ihkv_0}$, $u=e^{-hs}$,
and identify $s$ with the Laplace variable (with $k=\kappa/(|v_0|h) \in (\,-\pi/(|v_0|h), \pi/(|v_0|h)\,) \to (-\infty,\infty)$, \, $v_0={\tilde \sigma}_0|v_0|$). This yields the Fourier-Laplace transform\footnote{With Fourier-Laplace inverses $\chi_{1}(t,x)=\chi_{1}(t,x,v_0) = \xi_0 e^{-\xi_0t}\delta(x-v_0t)$ and $\chi_{2}(t,x)= \chi_{1}(t,x,-v_0)$.}
\beq
\label{C_limit}
 {\hat \chi}_1(s,k) =\lim_{h\to 0} {\bar \psi}_B[e^{-h(s+iv_0k)},h\xi_0] =\lim_{h\to 0} \frac{\xi_0he^{h(s+ikv_0)}}{e^{h(s+ikv_0)}-1+h\xi_0} =
\frac{\xi_0}{\xi_0+s+ikv_0}
\eeq
and ${\bar \psi}_B(u\zeta_2) \to {\hat \chi}_2(k,s) = \frac{\xi_0}{\xi_0+s-ikv_0}$.
The characteristic function (\ref{charfu_SRW}) converges with these scaling assumptions
to the Fourier-Laplace transform of the continuum-limit propagator
\beq
\label{Fourier_Laplace}
\begin{array}{clr}
\ds {\hat P}_c(k,s,v_0) & = \ds \lim_{h\to 0} h\, {\bar g}(e^{-sh},e^{-ikv_0h},e^{ikv_0h}) & \\ \\
 & = \ds 
\frac{1}{1-{\hat \chi}_{1}(k,s){\hat \chi}_{2}(k,s)}\left(\frac{1-{\hat \chi}_{1}(k,s)}{s+ikv_0}+ 
\frac{1-{\hat \chi}_{2}(k,s)}{s-ikv_0}
{\hat \chi}_{1}(k,s)\right) & \\ \\
 &  = \ds \frac{s+2\xi_0 -ik v_0}{s(s+2\xi_0)+k^2v_0^2} . &
 \end{array}
\eeq
The property ${\hat P}_c(0,s,v_0)=1/s$ shows the normalization of the propagator.
The drift term $\propto -ik$ generates bias and contains the Laplace transform of the expected position 
$ i\frac{\partial {\hat P}}{\partial k}(k,s)\big|_{k=0} = \frac{v_0}{s(s+2\xi_0)}$ 
which has the Laplace inverse 
\beq
\label{possijn}
\langle X_{v_0,t}^{(c)} \rangle = v_0\int_0^t \langle (-1)^{{\cal M}_P(\tau)} \rangle {\rm d}\tau = v_0 \int_0^t {\rm d}\tau \sum_{n=0}^{\infty}  e^{-\xi_0 \tau} \frac{(\xi_0 \tau)^n}{n!} (-1)^n =  \frac{v_0}{2\xi_0} (1-e^{-2\xi_0 t})
\eeq
where ${\cal M}_P(t) \in \mathbb{N}_0$ stands for the Poisson counting variable. Indeed (\ref{possijn}) is in agreement with the definition of the classical telegraph process  \cite{Kac1974,Orsingher1990,BogachevRatanov2011} defined by the random variable
$X_{v_0,t}^{(c)} = v_0\int_0^t (-1)^{{\cal M}_P(\tau)} {\rm d}\tau$, i.e. the velocity is reversed at the instants of Poisson events.
From (\ref{Fourier_Laplace}) we read off the partial differential equation governing in the space-Laplace domain
\beq
\label{tele_bias}
v_0\frac{d^2}{d x^2}{\hat P}_c(x,s,v_0)
-\frac{s(s+2\xi_0)}{v_0}{\hat P}_c(x,s,v_0) + \frac{s+2\xi_0}{v_0}\delta(x)-\frac{d}{dx}\delta(x) = 0
\eeq
where Dirac's $\delta$-distribution $\delta(x)=P_c(x,0,v_0)$ is the presumed initial condition. Eq. (\ref{tele_bias}) is a telegrapher's type equation with drift term $-\frac{d}{dx}\delta(x)$. 
Kac in his 1974 paper considers the symmetric 
combination (See \cite{Kac1974}, Eqs. (26), (30)) with the propagator
\beq
\label{walk_symmetric}
P_{tele}(x,t)= \frac{1}{2}\langle \delta(x-X^{(c)}_{v_0,t})
+\delta(x-X^{(c)}_{-v_0,t}) \rangle
= \frac{1}{2}[P_c(x,t,v_0)+P_c(x,t,-v_0)]
\eeq
which is an unbiased walk. The Fourier-Laplace 
transform then takes the form (canceling out the drift term) 
\beq
\label{drift_term}
{\hat P}_{tele}(k,s)= 
\frac{s+2\xi_0}{s(s+2\xi_0)+k^2v_0^2}
\eeq
which is the result reported by Kac \cite{Kac1974}. Worthy of mention is also the Gaussian limit $\xi_0 \to \infty$ while $\frac{v_0^2}{2\xi_0}=D$ is kept constant (`Kac limit') to standard Brownian motion (Wiener process) ${\hat P}_{tele}, {\hat P}_c \to 1/(s+Dk^2)$ removing the drift term in (\ref{Fourier_Laplace}). 
Expression (\ref{drift_term}) is an even function in $k$ and solves the classical telegrapher's (also called Cattaneo) equation (extensively studied in the literature \cite{BogachevRatanov2011}, and see the references therein) and writes in space-Laplace domain (\cite{Kac1974}, Eqs. (45), (46))
\beq
\label{telegraphers_eq}
v_0\frac{d^2}{d x^2}{\hat P}_{tele}(x,s)
-\frac{s(s+2\xi_0)}{v_0}{\hat P}_{tele}(x,s) + \frac{s+2\xi_0}{v_0}\delta(x) =0
\eeq
with initial condition 
$P_{tele}(x,t,v_0)\big|_{t=0}=\delta(x)$.
To depict the general structure of the propagators
consider the causal auxiliary Green's function $G(x,t)$ of the telegrapher's equation defined by
\beq
\label{telgreen}
\left[\frac{\partial^2}{\partial t ^2}+2\xi_0\frac{\partial}{\partial t} -
v_0^2\frac{\partial^2}{\partial x^2}\right]G(x,t)= \delta(x)\delta(t) 
\eeq
which has in Fourier-Laplace-space the representation $${\hat G}(k,s)=
\frac{1}{s^2+2s\xi_0+v_0^2k^2}= \frac{1}{\Lambda_1(k)-\Lambda_2(k)}\left(\frac{1}{s-\Lambda_1(k)}-\frac{1}{s-\Lambda_2(k)}\right)$$ where $\Lambda_{1,2}(k) = -\xi_0 \pm \sqrt{\xi_0^2-v_0^2k^2}$.
Laplace inversion yields in the $k$-$t$-space
\beq
\label{aux_greenf}
{\tilde G}(k,t)= e^{-\xi_0t}\frac{\sinh{\left(t\sqrt{\xi_0^2-v_0^2k^2}\right)}}{\sqrt{\xi_0^2-v_0^2k^2}}
\eeq
where we observe that $v_0{\tilde G}(k,t)\big|_{k=0} =
\langle X_{v_0,t}^{(c)} \rangle$ yields the expected position (\ref{possijn}).
Causality of $G(x,t)$ is ensured by $\xi_0>0$ where Eq. (\ref{telgreen}) in the $k$-$t$-space defines the response of a damped harmonic oscillator with spring constant $v_0^2k^2$ on an external forcing pulse.
The propagator (Fourier-Laplace inverse of (\ref{Fourier_Laplace})) has the representation
\beq
\label{occupationPDF}
P_c(x,t,v_0) = \left(\frac{\partial}{\partial t} - v_0\frac{\partial }{\partial x} + 2\xi_0\right)G(x,t) , \hspace{1cm} t>0
\eeq
with drift term 
$-v_0\frac{\partial G(x,t)}{\partial x}$ and $P_{tele}(x,t)= \left(\frac{\partial}{\partial t}+ 2\xi_0\right)G(x,t)$. 
We notice that auxiliary Green's function $G(x,t)$ is not a PDF, but the propagators
${\tilde P}_c(x,t,v_0), P_{tele}(x,t)$ are which is confirmed by
${\tilde P}_c(k,t,v_0)\big|_{k=0}={\tilde P}_{tele}(k,t)\big|_{k=0} =1$.
For further details we refer to the vast literature \cite{Goldstein1951,Kac1974,Orsingher1990,BeghinNiedduOrsingher2001} (and see also the references therein).
\subsection{Fractional Bernoulli SRW and its continuum limit}
\label{continuous-time_limits}
In this part we derive the space-time continuum limit for fat-tailed (FT) waiting densities taking us to a fractional generalization of the classical telegraph process. 
Consult \cite{MichelitschPolitoRiascos2021,TMM_FP_APR_fractiona_fract} for an outline of `well-scaled' continuum limits. Fat tailed waiting time density GFs have asymptotic representation (\ref{gen_feat}) with $\mu \in (0,1)$. More specifically we consider 
\beq
\label{ft-wait}
{\bar \psi}_{\mu}(u,\lambda)= \frac{u}{\lambda(1-u)^{\mu}+1} , \hspace{0.5cm} \mu \in (0,1) , \hspace{0.5cm} \lambda=\frac{1-p}{p} >0
\eeq
which refers to the ‘fractional Bernoulli process’ as a generalization of standard Bernoulli 
which is contained for $\mu=1$. The waiting time PDF which corresponds to
(\ref{ft-wait}) is the ‘discrete-time Mittag-Leffler distribution’ 
(of so-called ‘type B’) and a discrete version of the Mittag-Leffler density introduced recently \cite{PachonPolitoRicciuti2021}.
The fractional Bernoulli process indeed is a discrete-time version of the fractional Poisson point process
which was introduced by several authors
\cite{ScalasGorenfloMainardi2004,MainardiGorScalas2004,Laskin2003,GorenfloMainardi2011,BeghinOrsingher2009}.
The SRW associated to fractional Bernoulli (which we refer to as `Fractional Bernoulli SRW') yields with (\ref{mean_step_gen}) the GF of the average step
\beq
\label{en_fu_fract_ber}
{\bar \sigma}_{\mu,\lambda,{\tilde \sigma}_0}(u) = {\tilde \sigma}_0
\left(\frac{1 + \lambda(1-u)^{\mu-1}}{\lambda(1-u)^{\mu} +u+1}-1 \right).
\eeq
For $\mu=1$ this recovers relation (\ref{Bernoulli_genfu}) of the Bernoulli SRW. 
We can directly verify that the squirrel for
$t\to \infty$ escapes to infinity by the direction of ${\tilde \sigma}_0$ 
from the relation 
$\langle X_{\mu}(t)\rangle\big|_{t\to \infty}{\tilde \sigma}_0  ={\tilde \sigma}_0  {\bar \sigma}_{\mu}(u)\big|_{u\to 1} \to \,  \infty $ since $(1-u)^{\mu-1}$ is weakly singular at $u=1$ leading to the large time asymptotics with power-law escape (\ref{asymp_pos__large_t}).
Consider now the continuum limit yielding the expected velocity $\langle \, \sigma_t \rangle = \frac{d}{dt}\langle X_t\rangle$ in Laplace space
\beq
\label{scaled_frac_SRWlimit}
\begin{array}{clr}
\ds {\hat \sigma}(s) & = \ds
\lim_{h\to 0} h v_0 \left(\frac{1 + h^{-\mu}\xi_0^{-1}(1-e^{-hs})^{\mu-1}}{h^{-\mu}\xi_0^{-1}(1-e^{-hs})^{\mu} +e^{-hs}+1}-1 \right)  = v_0\frac{s^{\mu-1}}{s^{\mu}+2\xi_0} &
\end{array}
\eeq
where we have in (\ref{en_fu_fract_ber}) rescaled the constants in such a way that the limit $h\to 0$ exists, namely $\lambda(h)= h^ {-\mu}\xi_0^{-1}$, 
and the step size $|v_0|h$, with new constants $\xi_0>0$ (of units $[time^{-\mu}]$), and the directed velocity $v_0={\tilde \sigma}_0|v_0|$ independent of $h$.
The Laplace transform of the expected position then yields
\beq
\label{LT-expos}
{\hat X}^{(1)}(s) = \frac{{\hat \sigma}(s)}{s} = v_0\frac{s^{\mu-2}}{s^{\mu}+2\xi_0}.
\eeq
In this scaling limit, (\ref{ft-wait}) converges to the Laplace transform of the Mittag-Leffler density 
\beq
\label{ML_limit}
 \lim_{h\to 0}{\bar \psi}_{\mu}[e^{-hs},(\xi_0h^{\mu})^{-1}] = \lim_{h\to 0}
\frac{e^{-hs}}{ \xi_0^{-1} h^{-\mu} (1-e^{-hs})^{\mu}+1} = \frac{\xi_0}{\xi_0+s^{\mu}} .
\eeq
In the continuum limit the velocity is reversed at the instants of fractional Poisson events which we reconfirm subsequently. 
Now with above scaling assumptions we have ${\bar \psi}_{\mu}(u\zeta_{1,2},\lambda) \to \frac{\xi_0}{\xi_0+(s\pm ikv_0)^{\mu}}$ (and see (\ref{sum_up_yields_occupation_GF}), (\ref{charfu_SRW})) which yields for the Fourier-Laplace transform of the propagator
\beq
\label{Fourier_Laplace_Frac_Ber}
{\hat P}_{\mu}(k,s,v_0) = \lim_{h\to 0} h {\bar g}(e^{-hs},e^{-ikv_0},e^{ikv_0}) =
\frac{(s+ikv_0)^{\mu-1}[\xi_0+(s-ikv_0)^{\mu}]+\xi_0(s-ikv_0)^{\mu-1}}{[s^2+k^2v_0^2]^{\mu}+\xi_0(s+ikv_0)^{\mu}+\xi_0(s-ikv_0)^{\mu}}
\eeq
where necessarily ${\hat P}_{\mu}(0,s,v_0)=1/s$ (normalization of the propagator).
For $\mu=1$ this expression recovers (\ref{Fourier_Laplace}). We observe that (for real $s$) ${\hat P}_{\mu}(k,s,v_0) \in \mathbb{C}$ indicating a biased motion. This result shows that the continuum limit propagator solves a partial space-time fractional differential equation generalizing the `biased telegrapher's equation' (\ref{tele_bias}).
We introduce the Green's function of the `fractional telegrapher's equation'
\beq
\label{GF-grac-tele}
\left(\left[{\cal D}_{v_0}{\cal D}_{-v_0} \right]^{\mu}+\xi_0({\cal D}_{v_0}^{\mu}+ {\cal D}_{-v_0}^{\mu}) \right)G_{\mu}(x,t) =\delta(t)\delta(x)
\eeq
with ${\cal D}_{v_0} = \frac{\partial}{\partial t}+ v_0\frac{\partial}{\partial x}$.
The continuum limit propagator $P_{\mu}(x,t,v_0)$ then is represented by
\beq
\label{spatial-fracPDF}
P_{\mu}(x,t,v_0) = \left[{\cal D}_{v_0}^{\mu-1}(\xi_0+{\cal D}_{-v_0}^{\mu})+\xi_0{\cal D}_{-v_0}^{\mu-1}\right]G_{\mu}(x,t).
\eeq
The Green's function $G_{\mu}(x,t)$ and the resulting
propagator (\ref{spatial-fracPDF}) of the `fractional telegraph process' indeed merits further thorough analytical investigation. However, this is beyond the scope of the present paper. Instead we confine ourselves here to elaborate a few aspects and asymptotic features.
\\[1mm]
We confirm that
$i\frac{d}{dk}{\hat P}_{\mu}(k,s,v_0)\big|_{k=0}={\hat X}^{(1)}(s)= s^{2\mu-2}{\hat G}_{\mu}(k,s)\big|_{k=0}$ yields (\ref{LT-expos}) and recovers for $\mu=1$ the corresponding relation of the telegraph process of previous section. For $|s|$ small $ {\hat X}^{(1)}(s) \sim \frac{v_0}{2\xi_0} s^{\mu-2}$ the large time asymptotics for the FT case (\ref{asymp_pos__large_t}) is recovered (with rescaled constants)
\beq
\label{large-times}
\langle X_t \rangle \sim \frac{v_0}{2\xi_0} \frac{t^{1-\mu}}{\Gamma(2-\mu)} , \hspace{1cm} (t \to \infty).
\eeq
The squirrel escapes for $\mu \in (0,1)$ by a sublinear $t^{1-\mu}$-power law into the direction of $v_0$ (same direction as ${\tilde \sigma}_0$). For $\mu=1$ it recovers large time asymptotics $\frac{v_0}{2\xi_0}$ of the Poisson case (\ref{possijn}) of the classical telegraph process.
Laplace inversion of relation (\ref{scaled_frac_SRWlimit}) then yields 
\beq
\label{derivat}
\frac{d}{dt}\langle X_t\rangle = v_0 E_{\mu}(-2\xi_0 t^{\mu}) ,\hspace{1cm} t>0 
\eeq
where $E_{\mu}(z)$ indicates the standard Mittag-Leffler function defined by
(see e.g.\ \cite{MainardiGorScalas2004})
\beq
\label{Mittag_Leffler_function}
E_{\alpha}(z)= \sum_{m=0}^{\infty} \frac{z^m}{\Gamma(\alpha m+1)}.
\eeq
Therefore, we get (Laplace inversion of (\ref{LT-expos}))
\beq
\label{expect_pos_fracberoulli=}
\langle X_t\rangle = v_0 \int_0^t E_{\mu}(-2\xi_0 \tau^{\mu}){\rm d}\tau = v_0 t E_{\mu,2}(-2\xi_0 t^{\mu})
\eeq
fulfilling the initial condition $\langle X_0\rangle =0$
where the two-parameter Mittag-Leffler function
\beq
\label{two-param_ML}
E_{\alpha,\beta}(z)= \sum_{m=0}^{\infty} \frac{z^m}{\Gamma(\alpha m+\beta)}
\eeq
comes into play. Clearly relation (\ref{expect_pos_fracberoulli=})
is a fractional generalization of the telegraph process (\ref{possijn}) where the velocity directions change at arrival times of the fractional Poisson process which we show hereafter.
Taking into account the fractional Poisson state probabilities (probabilities of $n$ arrivals within $[0,t]$) \cite{Laskin2003}
\beq
\label{frac_Poiss}
\mathbb{P}[{\cal M}_{\mu}(t) = n] = \frac{(\xi_0 t^{\mu})^n}{n!}\frac{d^{n}}{dy^n} E_{\mu}(y)\big|_{y = -\xi_0 t^{\mu}}
\eeq
where ${\cal M}_{\mu}(t) \in \mathbb{N}_0$ stands for the fractional Poisson counting variable.
Thus $\langle (-1)^{{\cal M}_{\mu}(t)}\rangle = \sum_{n=0}^{\infty}\mathbb{P}[{\cal M}_{\mu}(t) = n](-1)^n =E_{\mu}(-2\xi_0t^{\mu})$ yields the expected position (\ref{expect_pos_fracberoulli=}) 
\beq
\label{expect_recover}
\langle X_t\rangle  = \langle \int_0^t (-1)^{{\cal M}_{\mu}(\tau)}{\rm d\tau}  \rangle = v_0 \int_0^t E_{\mu}(-2\xi_0\tau^{\mu}){\rm d}\tau
\eeq
with initial condition $\langle X_t\rangle\big|_{t=0}=0$ and $\frac{d}{dt}\langle X_t\rangle\big|_{t=0}=v_0$. For $\mu=1$ all expressions turn into the Poisson counterparts of the classical telegraph process.
\section{Anomalous diffusion}
\label{anomalous}
\subsection{The aging effect}
\label{aging_renewalp}
In this section we investigate the anomalous diffusive features
of the SRW for an arbitrary non-Markovian discrete-time renewal process 
${\cal N}(t)$ ($t \in \mathbb{N}_0$) defined in (\ref{discrete-time}). 
Recall the definition of anomalous diffusion \cite{MetzlerKlafter2001}: in a wide range of systems the mean square displacement scales with
a power law $\langle\, X_t^2\, \rangle \sim D_{\beta}t^{\beta}$ (where the displacement refers to the initial position with $X_0=0$)
and $D_{\beta}$ indicates the generalized diffusion coefficient (having units $[length^2\,time^{-\beta}]$). Then for 
$0<\beta<1$ the motion is subdiffusive, for $\beta=1$ normal diffusive and $\beta>1$ refers to superdiffusion where $\beta=2$ corresponds to
ballistic superdiffusion and $\beta>2$ to hyperballistic superdiffusion.
We especially focus on the variance
\beq
\label{variance_SQR}
{\cal V}(t) = \langle \left(X_t-\langle X_t 
\rangle\right)^2 \rangle = \langle X_t^2 \rangle - \langle X_t 
\rangle^2
\eeq
with the mean square displacement
\beq
\label{msd_discrte-time}
\langle X_t^2 \rangle = \sum_{k=1}^t\sum_{r=1}^t \langle (-1)^{{\cal N}(k)-{\cal N}(r)} \rangle = 
-t + 2K(t) 
\eeq
where we introduced the auxiliary quantity
\beq
\label{Kt}
K(t) = \sum_{r=1}^{t} 
\sum_{k=0}^{t-r} \langle (-1)^{{\cal N}_r(k)}\rangle .
\eeq
In this expression the integer counting variable ${\cal N}_r(\tau) \in \mathbb{N}_0$ appears and is defined by
\beq
\label{agin-discrete}
{\cal N}_{\tau}(t) = {\cal N}(t+\tau)-{\cal N}(\tau) ,\hspace{1cm} t,\tau = 0,1,2,,\ldots
\eeq
with initial condition ${\cal N}_{\tau}(t)\big|_{t=0}=0$. 
The quantity in (\ref{agin-discrete}) is the discrete-time version of the so called `aging renewal process' and is different from the original renewal process ${\cal N}(t)$ of (\ref{discrete-time}) if the latter is non-Markovian. For continuous times the aging renewal process was to our knowledge first introduced in \cite{GodrecheLuck2001} and for CTRW models based on aging renewal theory we refer to the references \cite{Barkai2003,BarkaiCheng2003,SchulzBarkaiMetzler2014}. 
We refer the counting process (\ref{agin-discrete}) to as `{\it discrete-time aging renewal process}' (DTARP) and call the (integer) variable $\tau$ `aging parameter'. Clearly ${\cal N}_{\tau=0}(t)={\cal N}(t)$ recovers the original renewal process. 
Intuitively, we infer that the events ${\cal N}_{\tau}(t) >1$ are drawn 
 from waiting-time density $\psi_t$ of the original renewal process ${\cal N}(t)$, however the density of the first event is different from $\psi_t$ and modifies (in the general non-Markovian case) the statistics.
We invite the reader to consult Appendices \ref{aging_discrete-time_renewal}, \ref{pertinent_DTARP} for detailed derivations and discussions of pertinent DTARP distributions and the related GFs which we employ extensively in the following. 
\subsection{Sibuya SRW}
\begin{figure*}[t]
\centerline{\includegraphics[width=0.6\textwidth]{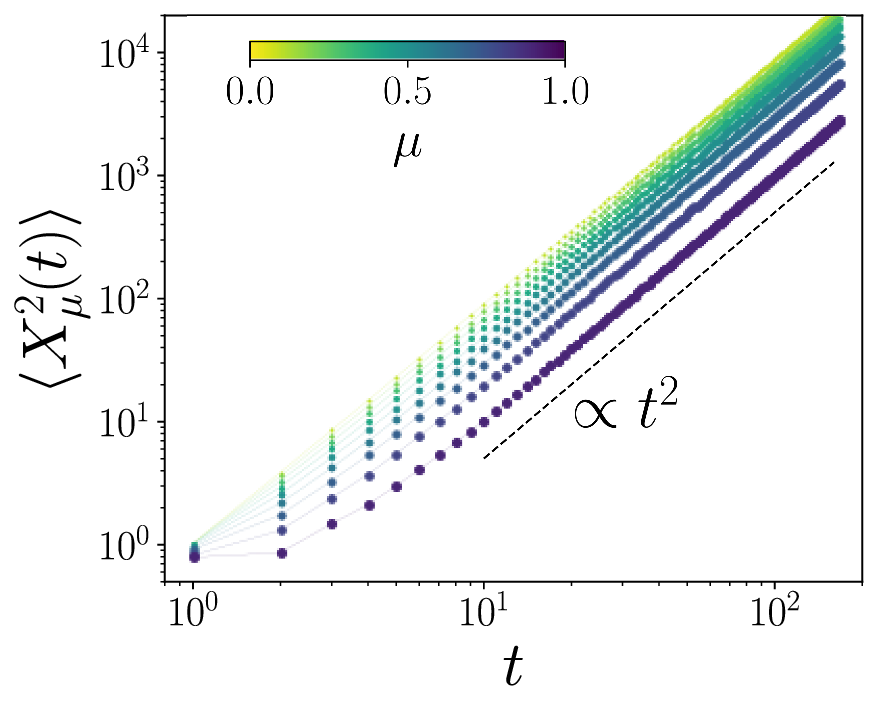}}
\vspace{-5mm}
\caption{\label{fig_x2mu} Numerical evaluation of $\langle X_{\mu}^2(t)\rangle$ of Eq.\ (\ref{MSD_sib}) as a function of $t$ for $\mu=0.1,0.2,\ldots,0.9$ codified in the colorbar. The dashed line represents the asymptotic power-law scaling  $\langle X_{\mu}^2(t)\rangle\propto t^2$.}
\end{figure*}

As a prototypical example of a non-Markovian SRW with strong aging effect we explore here the diffusive features of the `Sibuya SRW', i.e. the walk where the waiting times between the step reversals follow the Sibuya distribution. The Sibuya PDF has the GF
\beq
\label{Sibuya}
{\bar \psi}_{\mu}(u) = 1-(1-u)^{\mu} , \hspace{1cm} \mu \in (0,1).
\eeq
The Sibuya waiting-time PDF
has the form \cite{PachonPolitoRicciuti2021}
\beq
\label{Sibuya_PDF}
\psi_{\mu}(t) =  \frac{(-1)^{t-1}}{t!}\mu(\mu-1)\ldots (\mu-t+1) =  \frac{\mu\Gamma(t-\mu)}{\Gamma(1-\mu)\Gamma(t+1)}
\eeq
and is fat-tailed (FT), i.e. the expected waiting time
is infinite, $\frac{d}{du}{\bar \psi}_{\mu}(u)\big|_{u=1} \to \infty $, since
$\psi_{\mu}(t) \sim \mu t^{-\mu-1}/\Gamma(1-\mu)$ ($t\to \infty$). 
Using Eqs. 
(\ref{genfu}) with (\ref{GF_msd_general}) the GF of auxiliary quantity (\ref{Kt}) yields
\beq
\label{scaling_limit_Sibuya}
\begin{array}{clr}
\ds {\bar K}_{\mu}(u) & = \ds  (1-u)^{-3}\left(1-\frac{u\mu}{1-\frac{1}{2}(1-u)^{\mu}}\right)-\frac{1}{2}\frac{(1-u)^{\mu-2}}{1-\frac{1}{2}(1-u)^{\mu}}  &  \\ \\
 & = \ds (1-u)^{-3} -2\mu \, {\bar p}_{\mu,\mu+2}^1(2,u)+2\mu \, {\bar p}_{\mu,\mu+3}^1(2,u)+ {\bar p}_{\mu,2}^1(2,u)  & \\ \\
 & \sim \ds (1-\mu)(1-u)^{-3} +o[(1-u)^{-3}]   & (u \to 1-) 
 \end{array}
\eeq
where we introduced the GF of the discrete-time Prabhakar kernel \cite{MichelitschPolitoRiascos2021,TMM_FP_APR_Prab_Warsaw2021} (see Appendix \ref{Discrete_Prabhakar} for some details):
\beq
\label{PrabGF}
{\bar p}_{\mu,\nu}^{\gamma}(\lambda,u) = 
\frac{(1-u)^{-\nu}}{(1-\lambda(1-u)^{-\mu})^{\gamma}}.
\eeq
The continuous-time version of the Prabhakar kernel was first introduced by Giusti 
\cite{Giusti2020,Giusti-et-al-Prabhakar2020} (and see the references therein). 
Representation (\ref{Prabhakar_funct}) of the Prabhakar kernel allows us to invert (\ref{scaling_limit_Sibuya}) to arrive at the exact formula
(we employ the notation $p_{\mu,\nu}(\lambda,t)= p_{\mu,\nu}^1(\lambda,t)$)
\beq
\label{Kmu_exact}
K_{\mu}(t) =  
\frac{(t+1)(t+2)}{2}  +2\mu \, p_{\mu,\mu+3}(2,t) -2\mu \, p_{\mu,\mu+2}(2,t) + p_{\mu,2}(2,t).
\eeq
The mean square displacement then yields
\beq
\label{MSD_sib}
\langle X_{\mu}^2(t)\rangle  = 2K_{\mu}(t)-t = (t+1)^ 2+1
+4\mu \, p_{\mu,\mu+3}(2,t) -4\mu \, p_{\mu,\mu+2}(2,t) + 2p_{\mu,2}(2,t)
\eeq
where we verify that $K_{\mu}(0)=0$
and
$\langle X_{\mu}^2(0)\rangle=0$ as
$p_{\mu,\nu}(2,0)=-1$ and with $p_{\mu,\nu}(2,1)=2\mu-\nu$ we further confirm that necessarily $\langle X_{\mu}^2(1)\rangle=1$ (Appendix \ref{Discrete_Prabhakar}). 
Then to compute the variance we need the expected position
which we obtain as (see Eq. (\ref{gen_fu_of_position}))
\beq
\label{GF_Expect_Sibuya}
\langle X_{\mu}(t) \rangle = -{\tilde \sigma}_0(p_{\mu,2}(2,t)+1)
\eeq
with initial condition $\langle X_{\mu}(0) \rangle =0$. 
The variance of the Sibuya SRW then writes
\beq
\label{exact_Sibuya_variance}
{\cal V}_{\mu}(t) =  \ds 
\langle [X_{\mu}(t)]^2 \rangle - [\langle X_{\mu}(t)\rangle]^2  =(t+1)^2
+4\mu \, p_{\mu,\mu+3}(2,t) -4\mu \, p_{\mu,\mu+2}(2,t) - [p_{\mu,2}(2,t)]^2 
\eeq
which is an exact formula
where necessarily ${\cal V}_{\mu}(0)=0$.
For the large-time asymptotics this yields (see (\ref{as_t_small_and_large}) with (\ref{asymp_Prabhakar}))
\beq
\label{large_time}
K_{\mu}(t) \sim  \frac{(1-\mu)}{2}t^2 ,\hspace{1cm} (t\to \infty)
\eeq
and therefore
\beq
\label{MSD_asym}
\langle X_{\mu}^2(t)\rangle  = 2K_{\mu}(t)-t \sim (1-\mu) t^2 ,\hspace{1cm} (t \to \infty)
\eeq
which corresponds to superdiffusive ballistic $t^{2}$-scaling with generalized diffusion coefficient $D_{\mu}=1-\mu$
decreasing with increasing $\mu$ (i.e. for shorter waiting times between step reversals)
where this holds for the FT range $\mu \in (0,1)$. Such a ballistic large time asymptotics was also reported for fractional (continuous time) Cattaneo transport \cite{CompteMetzler1999}.
For large $t$ (\ref{GF_Expect_Sibuya}) has the asymptotics 
\beq
\label{large_t_exp_pos}
\langle X_{\mu}(t) \rangle \sim  {\tilde \sigma}_0 
\frac{t^{1-\mu}}{2\Gamma(2-\mu)} , \hspace{1cm} (t\to \infty)
\eeq
which is in agreement with Eq. (\ref{asymp_pos__large_t}) (with $A_{\mu}=1$ for Sibuya). We have then 
\beq
\label{Sib_meansquare}
\langle X_{\mu}(t)\rangle^2 \sim \frac{1}{4[\Gamma(2-\mu)]^2} t^{2-2\mu} \ll
t^{2} , \hspace{1cm} (t\to \infty).
\eeq
The large time asymptotics of the variance is therefore dominated by the mean square displacement (\ref{MSD_asym}), namely
\beq
\label{Sibuya_variance}
\begin{array}{clr} 
\ds {\cal V}_{\mu}(t) & \sim \langle X_{\mu}^2(t)\rangle \sim (1-\mu) t^2 , & \hspace{1cm} (t\to \infty) .
\end{array}
\eeq 
In the large-time limit the Sibuya SRW is superdiffusive with a ballistic $t^{2}$-law. This also holds true for the entire class of SRWs with fat-tailed waiting-time densities.
The ballistic scaling can be seen in Figure 
\ref{fig_x2mu} where we plot
$\langle X_{\mu}^2(t)\rangle$.
\\[1mm]
For $\mu=1$ the variance (\ref{Sibuya_variance}) is null
where this limit corresponds to the trivial deterministic counting process ${\cal N}_1(t)=t$ where the squirrel is trapped close to the departure site (this limit coincides with the limit $p=1$ previously discussed for the Bernoulli SRW -- see Section \ref{Bernoulli}).
The oscillating behavior can be extracted as the Sibuya GF then collapses to ${\bar \psi}_1(u)=u$ coinciding with Bernoulli for $p=1$. Then (\ref{scaling_limit_Sibuya}) yields
\beq
\label{mu_one}
{\bar K}_1(u) =\frac{u}{(1-u)^2(1+u)}
\eeq
in agreement with (\ref{Bernoulli_K}) for $p=1$.
Thus ${\bar X}_{\mu=1}^{(2)}(u)= u/[(1-u)(1+u)]$ which yields 
\beq
\langle X_{1}^2(t)\rangle = \sum_{r=0}^t \Theta(t-1-r)(-1)^r  =  \frac{1}{2}(1-(-1)^{t}) =
-{\tilde \sigma}_0\langle X_t\rangle 
\eeq
and is in agreement
with (\ref{meanquare_bern}). Therefore, we have indeed
${\cal V}_{\mu=1}(t)={\cal V}_{p=1}(t)_{Ber}=0$.
\\[1mm]
Further important is the `frozen limit' of infinite waiting times $\mu \to 0+$ which we already discussed in Section
\ref{large_times}.
\subsection{Anomalous diffusion of SRWs with broad and narrow waiting time densities}
\label{broad_densities}
To complement this part consider now the large time behavior for $\mu=1$ with expansion (\ref{infvar}). In this case, relation (\ref{genfu}) has the asymptotic expansion ($u\to 1-$)
\beq
\label{expanK}
{\bar K}(u) 
= \left\{\begin{array}{ll} \ds  \frac{1}{(1-u)^2}\left(1-\frac{A_1}{2}\right) + \frac{(\lambda-1)B_{\lambda}}{A_1}(1-u)^{\lambda-4} + o[(1-u)^{-2}], & 1<\lambda < 2 \\ \\ \ds
\frac{1}{(1-u)^2}\left[1+\frac{B_2}{A_1}-\frac{A_1}{2}\right]+o\left[(1-u)^{-2}\right], & \lambda =2 \end{array}\right.
\eeq
where the second line is consistent with the Bernoulli case (see Eq. (\ref{Bernoulli_K}) for $u\to 1-$). This takes us to the large time behavior
\beq
\label{msd_ano}
\langle X^2_t \rangle  = 2K(t)-t  \sim \left\{\begin{array}{ll} \ds (1-A_1)t  +
 \frac{2(\lambda-1)B_{\lambda}}{A_1} \frac{t^{3-\lambda}}{\Gamma(4-\lambda)}  \to  \frac{2(\lambda-1)B_{\lambda}}{A_1} \frac{t^{3-\lambda}}{\Gamma(4-\lambda)}\, \, ,  & \hspace{0.5cm} 1<\lambda < 2 \\ \\
\ds \frac{2B_2+A_1-A_1^2}{A_1} \, t  \, \, ,
 & \hspace{0.5cm} \lambda = 2 \end{array}\right.
\eeq
with the variance $2B_2+A_1-A_1^2= \sum_{k=1}^{\infty}[k-A_1]^2\psi(k)$.
For Bernoulli the last line yields $\langle X^2_t \rangle \sim t(1-p)/p$
in agreement with our previous results (see (\ref{meanquare_bern})).
Generally, for $\lambda=2$ this becomes a linear relation corresponding to normal diffusion whereas for broad waiting time densities ($\lambda \in (1,2)$) this is a superdiffusive law with a scaling exponent $1< 3-\lambda <2$.
\section{The SRW time-changed with a renewal process}
\label{examples}
Here we introduce a class of continuous-time walks by time-changing the SRW with an independent renewal process ${\cal M}(t) \in \mathbb{N}_0, t \in \mathbb{R}^+$ (i.e. an independent continuous-time counting process with IID interarrival times).  
This defines a biased continuous-time random walk which we call "continuous-time squirrel random walk" (CTSRW). 
The position of the squirrel in a CTSRW can be represented by the random variable $Y(t) \in \mathbb{Z}$ such that
\beq
\label{represent_thiswalk}
Y(t) = X_{{\cal M}(t)} =  X_{{\cal M}(t)-1} + \sigma_{{\cal M}(t)} , \hspace{1cm} t \in \mathbb{R}^{+} 
\eeq
with initial condition $Y(t)\big|_{t=0}=X_0=0$ and the increment
$$
\sigma_{{\cal M}(t)} ={\tilde \sigma}_0 \left[
(-1)^{{\cal N}[{\cal M}(t)]} -\delta_{{\cal M}(t),0} \right]
$$
where $X_{m\in \mathbb{N}_0}$ is the SRW (\ref{random_var}). 
In the CTSRW, the directed unit steps on $\mathbb{Z}$ are performed only
at the arrival time instants of the point process ${\cal M}(t) \in \mathbb{N}_0$ ($t \in \mathbb{R}^+$) defining a random clock (the operational time of the walk)
and $t$ is the (continuous) chronological time. 
In this time-change construction the step directions are reversed at the arrival times of the composed process ${\cal N}[{\cal M}(t)] \in \mathbb{N}_0$ which also is a point process defined on $t\in \mathbb{R}^+$.
We will see that the CTSRW is a different class to the class of continuous-time walks emerging by the continuum limits in the SRW, as considered in Section \ref{Diffusive_telepgraph}. On the contrary to the latter where the squirrel never waits, in the CTSRW the squirrel does not move during the inter-arrival time intervals of the point process ${\cal M}(t)$.
Consider now a PDF $f(r)$ supported on integers $r=0,1,2,\ldots$ and its GF ${\bar f}(u) =\sum_{r=0}^{\infty}f(r)u^r$ ($|u| \in [0,1]$). Its time-changed counterpart $f[{\cal M}(t)]$ has the mean
\beq
\label{expect_function}
F(t) = \langle f[{\cal M}(t)]\rangle= \sum_{m=0}^{\infty} \mathbb{P}[{\cal M}(t)=m] f(m) , \hspace{1cm} t \in \mathbb{R}^{+}
\eeq
which depends on the continuous time $t$.
In this relation, the state probabilities $P_m(t)=\mathbb{P}[{\cal M}(t)=m] $ (probabilities for $m=0,1,2,\ldots$ arrivals
within the continuous time interval $[0,t]$) come into play.
Let ${\hat \eta}(s)$ denote the Laplace transform of the interarrival time density $\eta(t)$ of the point process ${\cal M}(t)$. Then, from conditioning arguments we have
\beq
\label{state_prob_LT}
{\hat P}_m(s) = \int_0^{\infty} e^{-st} \mathbb{P}[{\cal M}(t)=m] {\rm d}t = 
\frac{1-{\hat \eta}(s)}{s} [{\hat \eta}(s)]^m.
\eeq
Thus, the Laplace transform of (\ref{expect_function}) yields
\beq
\label{LT_expect_function}
{\hat F}(s) = \frac{1-{\hat \eta}(s)}{s}\sum_{m=0}^{\infty} f(m) [{\hat \eta}(s)]^m = \frac{1-{\hat \eta}(s)}{s}{\bar f}[{\hat \eta}(s)]
\eeq
where ${\bar f}[{\hat \eta}(s)]$ is the GF with argument $u={\hat \eta}(s)$.
Eq. (\ref{LT_expect_function}) relates the GFs of functions defined on integer times with the Laplace transforms of their time-changed means. By using this general result we can represent the propagator of the CTSRW
(see (\ref{sum_up_yields_occupation_GF}), (\ref{charfu_SRW})) in Fourier-Laplace space as
\beq
\label{CTSRW_propagator}
 {\hat Q}(\kappa,s) = \frac{1-{\tilde \eta}(s)}{s}
 {\bar g}[{\tilde \eta}(s),e^{-i\kappa{\tilde \sigma_0}},e^{i\kappa{\tilde \sigma_0}}] , \hspace{1cm} \kappa \in (-\pi,\pi)
\eeq
where normalization of the propagator is confirmed ${\hat Q}(0,s)=\frac{1}{s}$
by using 
${\bar \psi}[{\tilde \eta}(s)]\big|_{s=0}= {\bar \psi}(1)=1$. Expression (\ref{CTSRW_propagator}) may serve as a point of departure for a wide field of applications of the CTSRW model.
\\[1mm]
As a useful example for the following consider the Laplace transform of the mean of the time-changed discrete Prabhakar kernel (see (\ref{PrabGF_pp}))
\beq
\label{time_changed_Ptab}
\int_0^{\infty}e^{-st}\langle
p_{\mu,\nu}(\lambda, {\cal M}(t) \rangle {\rm d}t = \frac{1-{\hat \eta}(s)}{s}{\bar p}_{\mu,\nu}(\lambda,{\hat \eta}(s))= \frac{1}{s}\,\frac{[1-{\hat \eta}(s)]^{1-\nu}}{1-\lambda[1-{\hat \eta}(s)]^{-\mu}} = \frac{1}{s}\, {\bar p}_{\mu,\nu-1}[\lambda,{\hat \eta}(s)].
\eeq
We can see here directly the long-time asymptotics. Consider a FT density
with expansion (\ref{small-s}). Then (\ref{time_changed_Ptab}) behaves
for $s\to 0$ as 
\beq
\label{as-LT}
\frac{1}{s}\, {\bar p}_{\mu,\nu-1}[\lambda,{\hat \eta}(s)] \sim -\frac{b^{1+\mu-\nu}}{\lambda}
s^{-\alpha(\nu-\mu-1)-1} 
\eeq
and therefore
\beq
\label{larg-timeprab}
\langle\, 
p_{\mu,\nu}[\lambda, {\cal M}(t)] \, \rangle \sim -\frac{b^{1+\mu-\nu}}{\lambda} \frac{t^{\alpha(\nu-\mu-1)}}{\Gamma(\alpha[\nu-\mu-1]+1)}, \hspace{1cm} (t \to \infty).
\eeq
This expression is the time changed version of the large time 
asymptotics (\ref{as_t_small_and_large}) and both coincide for LT case $\alpha=1$ (when $b=1$). For $\alpha \in (0,1)$ the $t^{\alpha(\nu-\mu-1)}$-Prabhakar power-law is slowing down the ballistic diffusive Sibuya SRW scaling (\ref{Sibuya_variance}), i.e. when we subordinate the Sibuya SRW to a
renewal process with fat-tailed waiting time density (see (\ref{small-s})).
This slowdown is a consequence of the long waiting intervals where 
the squirrel does not move.
We consider this issue in more detail subsequently (relation (\ref{LT_MSD_timechanged})).
\\[2mm]
With these remarks we can write for 
the expected CTSRW position 
\beq
\label{exp_pos}
\ds \langle Y(t) \rangle =  \sum_{m=0}^{\infty} \mathbb{P}[{\cal M}(t)=m]\langle X_{m}\rangle  , \hspace{1cm} t \in \mathbb{R}^{+}
\eeq
which has then the Laplace transform (see (\ref{gen_fu_of_position}))
\beq
\label{expect_pos_Laplace}
\begin{array}{clr}
\ds {\hat Y}^{(1)}(s) & = \ds
\frac{1-{\hat \eta}(s)}{s} {\bar X}^{(1)}[{\hat \eta}(s)] = \frac{1}{s}{\bar \sigma}[{\hat \eta}(s)] & \\ \\
  & = \ds\frac{{\tilde \sigma}_0}{s} \left(\frac{1-{\bar \psi}[{\hat \eta}(s)]}{[1-{\hat \eta}(s)](1+{\bar \psi}[{\hat \eta}(s)])} -1 \right). &
  \end{array}
\eeq
We point out that 
${\bar \psi}[{\hat \eta}(s)]$ is the Laplace transform of the waiting-time density of the composed counting process ${\cal N}[{\cal M}(t)]$ (see \cite{ADTRW2021,OrsingherPolito2012a,OrsingherPolito2012b} for details).
In order to explore how the time change affects anomalous diffusion we consider subsequently the
Laplace transform ${\hat Y}^{(2)}(s)$ of the CTSRW mean square displacement which takes the form (see (\ref{LT_expect_function}) and 
(\ref{genfu}) with (\ref{GF_msd_general}))
\beq
\label{MDS_time-changed}
{\hat Y}^{(2)}(s) = \frac{1}{s}
\left(2[1-{\hat \eta}(s)]{\bar K}[{\hat \eta}(s)] - 
\frac{{\hat \eta}(s)}{1-{\hat \eta}(s)}\right).
\eeq
\subsection{Bernoulli SRW time-changed with an arbitrary renewal process}
\label{Bernoulli_CTSRW}
As an example consider the Bernoulli SRW subordinated to an independent arbitrary renewal process. Using (\ref{expect_pos_Laplace}) with (\ref{Bernoulli_genfu}) we can write the Laplace transform of the expected position as
\beq
\label{expect_pos}
{\hat Y}_p^{(1)}(s) = {\tilde \sigma}_0\frac{1-2p}{2p s}\, \frac{2p{\hat \eta}(s)}{[1-{\hat \eta}(s)(1-2p)]}.
\eeq
Note that the part ${\hat g}_p(s)= \frac{2p{\tilde \eta}(s)}{[1-{\hat \eta}(s)(1-2p)]}$ is the Laplace transform of a density $g_p(t)$ as ${\hat g}(s)\big|_{s=0}=1$
and in the expression (\ref{expect_pos}) the Laplace transform ${\hat g}_p(s)/s$
of its cumulative distribution $G_p(t) =\int_0^tg_p(\tau){\rm d}\tau$ is contained.
As a proto-typical example for long waiting times with fat-tailed waiting time density we consider the time fractional Poisson process with ${\hat \eta}_{\alpha,\xi}(s)=\frac{\xi}{\xi+s^{\alpha}}$ ($\xi>0$), thus the density $g_p(t)$
can then be identified with a Mittag-Leffler density and $G_p(t)= 1-E_{\alpha}(-2p\xi t^{\alpha})$ with the
Mittag-Leffler distribution.
The function $E_{\alpha}(z)$ stands for the standard Mittag-Leffler function (\ref{Mittag_Leffler_function}).
Notice that $E_{1}(z)=e^{z}$, which reflects the fact that for $\alpha=1$ all relations turn into the SRW time-changed with the standard Poisson.
The expected position takes the form
\beq
\label{expoect_pos_Bern}
\langle \, Y_p(t)\, \rangle  = {\tilde \sigma}_0\frac{1-2p}{2p} G_p(t)
\eeq
where for $t \to \infty$ we have $G_p(t) \to 1$ in agreement with the subsequent asymptotic relation (\ref{asymp_scaling}) for $\mu=1$ with $A_1=1/p$.
\subsection{Large time asymptotics of the CTSRW}
\label{lrage_time}
\begin{figure*}[t]
\centerline{\includegraphics[width=0.7\textwidth]{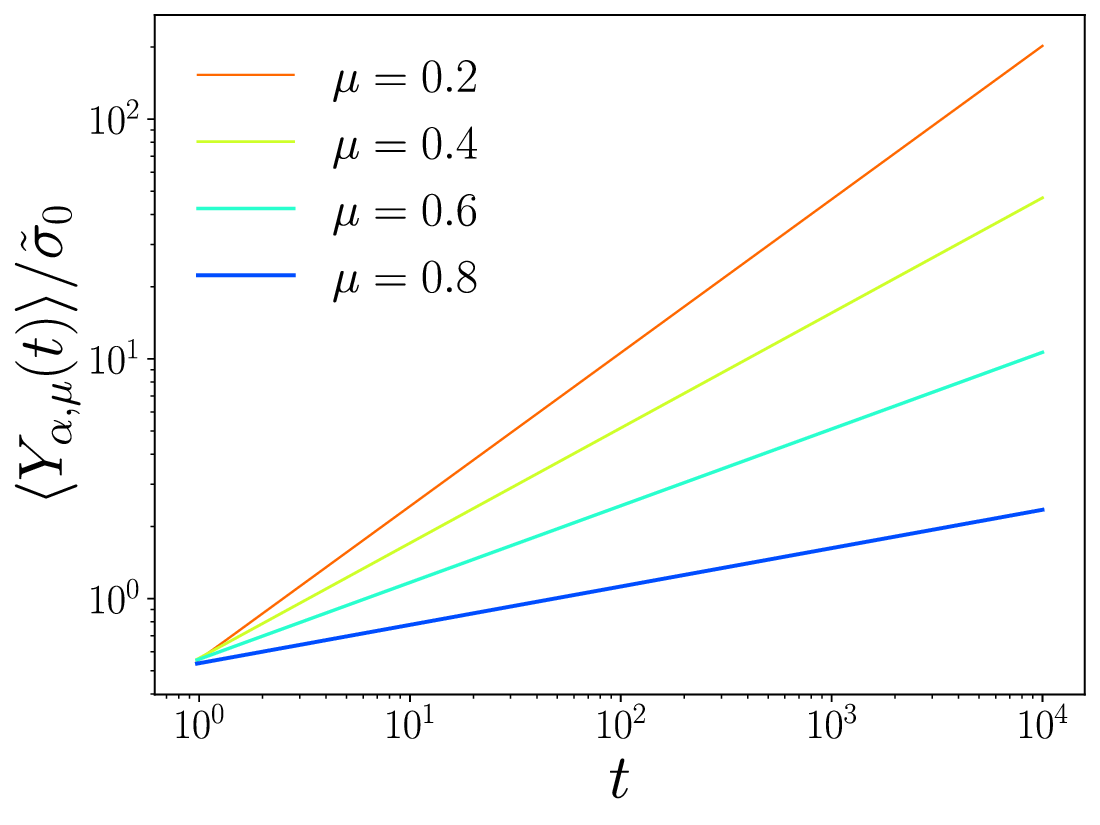}}
\caption{\label{fig1} Large time asymptotic scaling of $\langle Y_{\alpha,\mu}(t)\rangle$ in Eq. (\ref{asymp_scaling}) for the expected position with $\alpha=0.8$, $A_{\mu}=1$, $b=1$ for different values of $\mu$.}
\end{figure*}
\begin{figure*}[t]
\centerline{\includegraphics[width=0.75\textwidth]{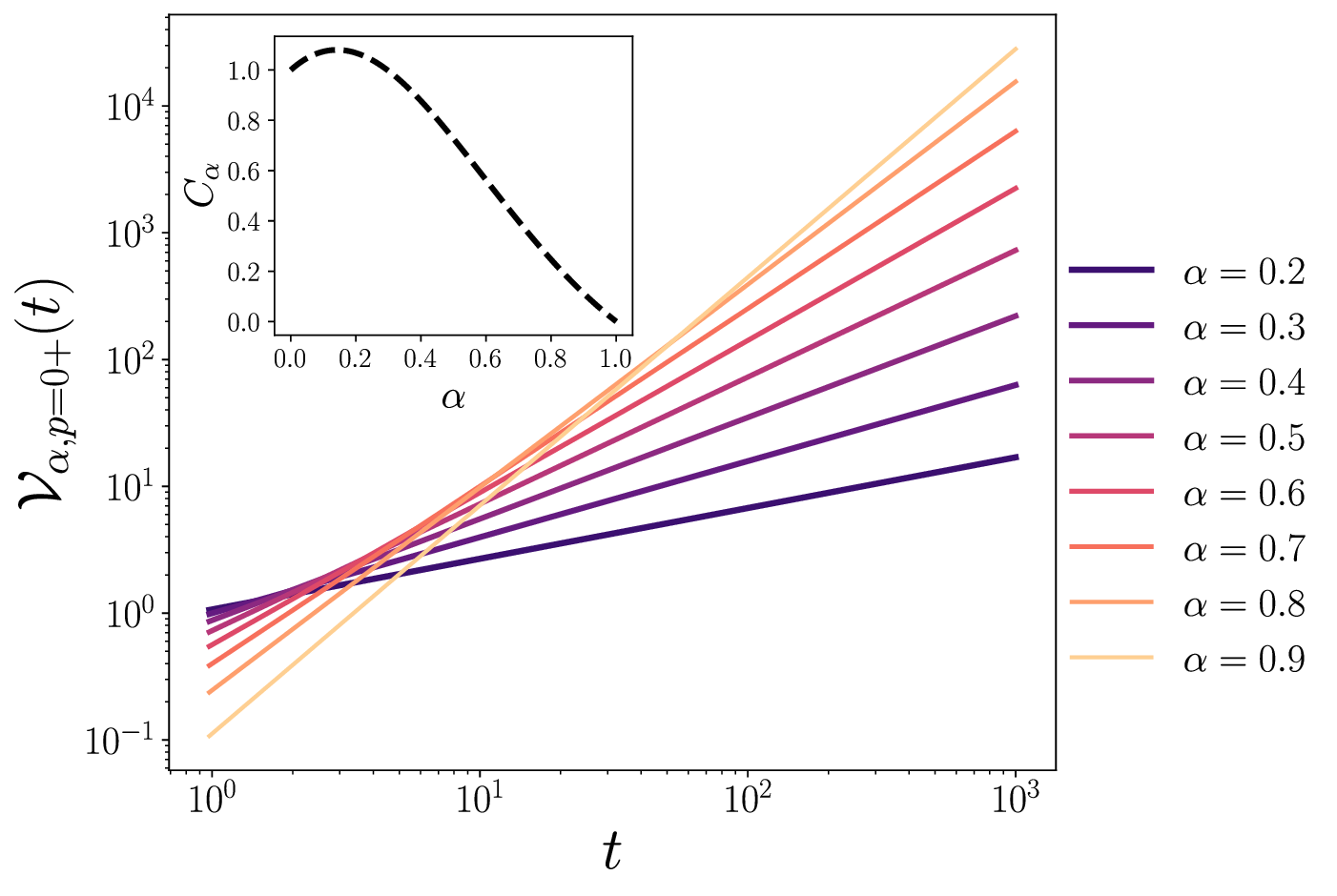}}
\caption{\label{fig3} Asymptotic limit of the variance ${\cal V}_{\alpha,p=0+}(t) $ in Eq. (\ref{mupzero_equal}) as a function of $t$ for different values of $\alpha$ and $b=1$. The inset shows $C_{\alpha}$ versus $\alpha$.}
\end{figure*}

For the large time behavior of the diffusive features we consider the Laplace transform of the waiting-time PDF which has the expansion
\beq
\label{small-s}
{\hat \eta}_{\alpha}(s) \sim 1- b s^{\alpha} +o(s^{\alpha}) , \hspace{1cm} \alpha \in (0,1], \hspace{0.5cm} (s\to 0)
\eeq
with $b>0$ and for $\alpha \in (0,1)$ the density $\eta_{\alpha}(t)$ is fat-tailed and for
$\alpha=1$ we confine us here to the case of narrow (light-tailed) waiting time densities $\eta(t)$ of the point process.
Expanding (\ref{expect_pos_Laplace}) for small $|s|$ by accounting for (\ref{genera_as_uto_one}) we arrive at
\beq
\label{Beroulli_gen}
\ds Y_{\alpha,\mu}^{(1)}(s) \sim \ds  \left\{\begin{array}{llr} \ds  \frac{{\tilde \sigma}_0}{2} A_{\mu}b^{\mu-1}s^{\alpha(\mu-1)-1} , & \hspace{1cm} \mu \in (0,1) , &  \\ \\
 \ds {\tilde \sigma}_0 \frac{(A_1-2)}{2s} -{\tilde \sigma}_0\frac{B_{\lambda}}{2}b^{\lambda-1}s^{\alpha(\lambda-1)-1} ,& \hspace{1cm} \mu=1,\hspace{1cm} \lambda \in (1,2]
\end{array} \right. \hspace{2cm} \alpha \in (0,1] .
\eeq 
By Laplace inversion we get for $t$ large
\beq
\label{asymp_scaling}
\ds \langle Y_{\alpha,\mu}(t)\rangle \sim  \ds  \left\{\begin{array}{llr}  \ds \frac{{\tilde \sigma}_0}{2}A_{\mu}b^{\mu-1} \frac{t^{\alpha(1-\mu)}}{\Gamma[1+\alpha(1-\mu)]}  , & \hspace{1cm} \ds \mu \in (0,1) , &  \\ \\
\ds   {\tilde \sigma}_0\frac{(A_1-2)}{2}  -{\tilde \sigma}_0\frac{B_{\lambda}}{2}b^{\lambda-1}\frac{t^{-\alpha(\lambda-1)}}{\Gamma(1-\alpha[\lambda-1])}  \to  {\tilde \sigma}_0\frac{(A_1-2)}{2} ,  & \hspace{1cm} \mu=1 ,& 1<\lambda \leq 2
\end{array} \right. 
\eeq 
where $A_{\mu}=1$ in the case of a time changed Sibuya SRW.
For $\mu=1$ the same asymptotic value ${\tilde \sigma}_0\frac{(A_1-2)}{2}$
as in the discrete-time case is approached (see (\ref{asymp_pos__large_t})) by a $t^{-\alpha(\lambda-1)}$ power law when $\alpha(\lambda-1) \neq 1$, and at least exponentially for $\lambda=2$, $\alpha=1$.
In contrast, when $\mu<1$, the squirrel escapes with a power law to the same direction as
${\tilde \sigma_0}$ whereas for $\mu=1$ it remains localized close to the departure site due to the oscillatory motion as in the discrete-time case (see (\ref{asymp_pos__large_t})).
We depict the large time asymptotic power-law behavior for four values of $\mu$ and $\alpha=0.8$ in Fig. \ref{fig1}. One can see that for increasing $\mu$ where the squirrel more often changes
the step directions the escape becomes slower.
On the other hand,
for small $\mu$ (long waiting times between the reversals of step directions) the squirrel escapes faster.
Now we finally consider (\ref{MDS_time-changed}) and without loss of generality the Sibuya CTSRW, i.e. $\mu \in (0,1)$ and $\alpha \in (0,1]$ ($A_{\mu}=1$).
Then we get
\beq
\label{LT_MSD_timechanged}
 {\cal V}_{\alpha,\mu}(t) \sim \langle Y_{\alpha,\mu}^2(t)\rangle \sim
 \frac{2(1-\mu)}{b^2}\frac{t^{2\alpha}}{\Gamma(2\alpha+1)} , \hspace{1cm} \mu\in (0,1),\hspace{0.3cm} \alpha \in (0,1] \hspace{0.5cm} (t\to \infty).
\eeq
We can obtain this result also directly by
using (\ref{larg-timeprab}) with (\ref{exact_Sibuya_variance}).
This relation is the time changed version of the ballistic Sibuya square law (\ref{Sibuya_variance}) and coincides with the large-time asymptotics of the mean square displacement. Taking into account (\ref{asymp_scaling}) we have $(\langle Y_{\alpha,\mu} \rangle)^2 \propto t^{2\alpha(1-\mu)} \ll  \langle Y_{\alpha,\mu}^2(t)\rangle$, thus the variance asymptotically is dominated by the mean square displacement. Be aware that the asymptotic formula (\ref{LT_MSD_timechanged}) is modified in the frozen limit for $\mu\to 0+$ considered at the end of this section.
We
identify two regimes of anomalous diffusion: For $0<\alpha< 1/2 $ the CTSRW is subdiffusive,
for $\alpha=1/2$ it is normal-diffusive, and for $1/2< \alpha \leq 1$ superdiffusive.
Be reminded that these regimes exist for $0<\mu<1$. For $\mu=1$ the large-time behavior is identical with the Bernoulli SRW subordinated to a renewal process (considered in the following).
\\[2mm]
{\it Bernoulli CTSRW} \\[2mm]
Considering (\ref{Ber_msdgenfu}) and (\ref{LT_expect_function}), we have 
the relation
\beq
\label{Bernoulli-time_change}
{\hat Y}_{\alpha}^{(2)}(s)_{Ber} = 
\frac{{\hat \eta}_{\alpha}(s)}{s[1-{\hat \eta}_{\alpha}(s)]}\frac{[1+{\hat \eta}_{\alpha}(s)(1-2p)]}{[1-{\hat \eta}_{\alpha}(s)(1-2p)]} \sim \frac{1-p}{pb}s^{-\alpha-1} , \hspace{1cm} (s \to 0)
\eeq
and therefore the mean square displacement scales as
$\langle Y_{\alpha}^2(t)_{Ber} \rangle \sim 
\frac{1-p}{pb} t^{\alpha}/\Gamma(\alpha+1)$. On the other hand the Laplace transform of the mean position
writes (see (\ref{later_use}) and (\ref{LT_expect_function}))
\beq
\label{LT_exptect-Ber}
{\hat Y}^{(1)}(s)_{Ber} \sim \frac{1-2p}{2p} \frac{1}{s} , \hspace{1cm} (s \to 0)
\eeq
thus the large time asymptotics
$\langle Y_{\alpha}(t)_{Ber} \rangle  \sim \frac{1-2p}{2p}$ is constant independent of $\alpha$ and identical with the  discrete-time case. Therefore
\beq
\label{subdiff_ber}
{\cal V}_{\alpha}(t)_{Ber} \sim 
\langle Y_{\alpha}^2(t)_{Ber} \rangle \sim \frac{1-p}{pb} \frac{t^{\alpha}}{\Gamma(\alpha+1)} ,\hspace{1cm} (t\to \infty)
\eeq
and is the time-changed version of the (linear) large-time asymptotics of (\ref{meanquare_bern}). We see that the effect of the time change is here
a subdiffusive power-law as consequence of long waiting times between the events of the point process. 
\\[2mm]
{\it Frozen limit} $p\to 0+$, $\mu\to 0+$\\[2mm]
Let us consider the  $p\to 0+$, $\mu\to 0+$, respectively (see again Section \ref{large_times}) where the squirrel is trapped for an infinitely long waiting time in the frozen regime. In the time
changed case this refers to the strict walk $Y(t)={\tilde \sigma}_0{\cal M}(t)$ without step reversals for finite $t$.
Therefore,
$\langle Y_{\alpha,p=0+}^2(t)\rangle = \langle Y_{\alpha,\mu=0+}^2(t)\rangle = \frac{2}{b^2}\frac{t^{2\alpha}}{\Gamma(2\alpha+1)}$ which is the time-changed version of (\ref{expand_p})). The variance then is
\beq
\label{mupzero_equal}
{\cal V}_{\alpha,p=0+}(t) = {\cal V}_{\alpha,\mu=0+}(t) = \langle \,[{\cal M}(t)]^2\, \rangle -\langle \,{\cal M}(t)\, \rangle^2\sim 
\frac{C_{\alpha}}{b^2}\, t^{2\alpha} , \hspace{1cm} (t\to \infty)
\eeq
with $C_{\alpha} = \frac{2}{\Gamma(2\alpha+1)}-\frac{1}{(\Gamma(\alpha+1))^2}$
where this relation holds for frozen limits of both Bernoulli and Sibuya CTSRW.
Contrary to the frozen limit in the SRW, the time changed version is not deterministic with a non-vanishing variance as its operational time is ${\cal M}(t)$.
In Fig. \ref{fig3}, we show the numerical values of the asymptotic limit of the variance ${\cal V}_{\alpha,p=0+}(t)$ in Eq. (\ref{mupzero_equal}) for different values of $\alpha$ showing the power-law scaling of the variance for large $t$. We also depict the $\alpha$-dependence of the multiplier $C_{\alpha}$ which is approaching zero for $\alpha \to 1$. This result reflects the vanishing variance of the deterministic SRW for $p, \mu \to0+$ which occurs here in the large-time limit.
\section{Conclusions}
We have introduced a semi-Markovian discrete-time generalization of the telegraph process, the `squirrel random walk'. Except for the Bernoulli SRW with geometric waiting times this walk is non-Markovian. We have derived an exact formula for the GF of the SRW characteristic function which determines uniquely the propagator (Eqs. (\ref{sum_up_yields_occupation_GF}), (\ref{charfu_SRW})).
We analyzed the diffusive limits to continuous space and time. For the Bernoulli SRW this leads to the standard telegraph process where the propagator is governed by a telegrapher's equation with drift.

We considered also a non-Markovian generalization, the fractional Bernoulli SRW. In this case the diffusive limit yields a fractional generalization
of the telegraph process where the propagator is governed by 
the `fractional telegrapher's equation' (see (\ref{GF-grac-tele}), (\ref{spatial-fracPDF})). This motion and propagator calls for further analytical investigation.
By taking into account the `aging effect' which naturally comes into play in non-Markovian SRWs, we explored the discrete-time counterpart ${\cal N}_{\tau}(t)$ of the so called `aging renewal process' (Eq. (\ref{agin_renewal})).
We derived explicit formulae for the mean square displacement and variance for the Sibuya SRW in terms of discrete-time Prabhakar kernels.
For the large-time limit we obtained a ballistic superdiffusive $t^2$-law 
as a hallmark of fat-tailed waiting time densities.

For SRWs governed by broad waiting-time PDFs (see (\ref{msd_ano})) we obtained superdiffusive large time scaling which turns into normal diffusion for narrow waiting time PDFs.
We also introduced time changed versions of the SRW (i.e. the SRW subordinated to independent point renewal processes) defining a class, the
`continuous-time squirrel random walk' and analyzed some large time asymptotic features. Pertinent candidates for future research include SRWs involving renewal processes being generalizations of fractional Poisson such as Prabhakar generalizations which have attracted recently a lot of interest \cite{CahoyPolito2013,MichelitschRiascos2020a,MichelitschRiascos2020b,MichelitschPolitoRiascos2021,TMM_FP_APR_Prab_Warsaw2021,Giusti-et-al-Prabhakar2020,Giusti2020}.
Generally the SRW model has a rich potential of applications and generalizations. 
For instance it may be generalized to multiple spatial dimensions or applied to random motions in graphs where at event instants of a renewal process a random choice is made among a certain set
of possible states or
positions (for models related to the latter, see 
\cite{RiascosBoyerHerringerMateos2020,SingGorskaSandev2022}). 
Finally, discrete time versions of aging renewal processes including variants with multivariate aging parameters and applications to random walks are interesting research directions as well.
\begin{appendix}
\section{Appendix}
\subsection{Discrete-time aging renewal process}
\label{aging_discrete-time_renewal}
For continuous-times the `aging renewal process' was introduced and analyzed in 
\cite{GodrecheLuck2001}
and for related aging continuous-time random walk models consult \cite{Barkai2003,BarkaiCheng2003,SchulzBarkaiMetzler2014}.
We introduce the `discrete-time aging renewal process (DTARP)' as
\beq
\label{agin_renewal}
N_{\tau}(t) = {\cal N}(t+\tau)-{\cal N}(\tau) , \hspace{1cm} t,\tau =\{0,1,2,\ldots\} \in \mathbb{N}_0
\eeq
where ${\cal N}(t)$ is the counting process (\ref{discrete-time}) and initial condition $N_{\tau}(t)\big|_{t=0} = 0$.
The counting process $N_{\tau}(t)$
represents the number of events occurring in the time interval $\{\tau+1,\ldots,\tau +t \}$ and for $\tau=0$ it recovers ${\cal N}_0(t)={\cal N}(t)$. 
Let $\{J_n\}$ be the renewal chain (\ref{renewal_chain}).
Following \cite{GodrecheLuck2001},
we define the `forward recurrence time' $E_{n,\tau}= \min(J_{n+1} \in \mathbb{N}_0 : J_{n+1}-\tau >0 )$ as the time interval from $\tau$ to the first event occurring later than $\tau$. We call the integer time $\tau$ `aging parameter'. 
Consider now for a fixed $n\in \mathbb{N}_0$ ($J_n\leq \tau)$ the probability that $E$ has a certain value 
$E_{n,\tau}=J_{n+1}-\tau= t$, i.e. the probability for
the first arrival occurring later than $\tau$ given that $n$ arrivals
are observed up to $\tau$ with $J_n \leq \tau < J_{n+1}$.
Denoting with $\mathbb{P}(A|B)$ the probability of $A$ conditional to $B$, this probability is defined as
\beq
\label{fE}
\begin{array}{clr}
f_E(\tau,t,n) &= \mathbb{P}[E_{n,\tau}=t|{\cal N}(\tau)=n] , &
 \tau,t,n \in \mathbb{N}_0\\ \\ & = \langle \, \Theta(J_n,\tau,J_{n+1}) \delta_{t,J_{n+1}-\tau} \, \rangle   &
\end{array}
\eeq
where $\delta_{i,j}$ denotes the Kronecker symbol and with the step function $\Theta(a,r,b)$ defined in
(\ref{gen_Theta}), i.e.
 $\Theta(J_n,\tau,J_{n+1})=1$ for $\tau \in [J_n,J_{n+1}-1]$, i.e. when ${\cal N}(\tau)=n$ and 
 $\Theta(J_n,\tau,J_{n+1})=0$
 elsewhere.
Then, the double GF of (\ref{fE}) is given by 
\beq
\label{doublegenfe}
\begin{array}{clr}
\ds {\bar f}_E(w,u,n) & = \ds \left\langle \, \sum_{r=0}^{\infty}w^r \Theta(J_n,r,J_{n+1}) \sum_{s=0}^{\infty} u^s\delta_{s,J_{n+1}-r} \,  \right\rangle  & n \in \mathbb{N}_0\\ \\
 & = \ds \left \langle u^{J_{n+1}} \sum_{r=J_n}^{J_{n+1}-1} u^{-r} w^r \right\rangle &  \\ \\
 & = \ds \left  \langle u^{J_{n+1}-J_n} w^{J_n} \sum_{r=0}^{J_{n+1}-J_n-1} w^r u^{-r}  \right\rangle   & \\ \\
 & = \ds \left  \langle w^{J_n} \right\rangle
 \frac{\langle u^{\Delta t_{n+1}}\rangle -\langle w^{\Delta t_{n+1}}\rangle}{1-\frac{w}{u}}  = [{\bar \psi}(w)]^n 
 \frac{u}{u-w}[{\bar \psi}(u)-{\bar \psi}(w)] &  
 \end{array}
\eeq
with $J_{n+1}-J_n=\Delta t_{n+1}$ (see (\ref{renewal_chain})) and where the IID feature of the $\Delta t_j$ with (\ref{gen_fu}) is used (and keep in mind the conjugations $u \leftrightarrow t$ and $w \leftrightarrow \tau$).
We observe that ${\bar f}_E(w,u,n)\big|_{u=0}=0= f_E(\tau,0,n)$ reflecting $N_{\tau}(0)=0$.
Further it is noteworthy that
\beq
\label{noteworthy}
{\bar f}_E(w,u,n)\big|_{u=1} = [{\bar \psi}(w)]^n \frac{1}{1-w}[1-{\bar \psi}(w)]
\eeq
is the GF of the state probabilities of the original renewal process 
$\sum_{t=0}^{\infty} f_E(\tau,t,n) =\langle \Theta(J_n,\tau,J_{n+1}) \rangle = \Phi^{(n)}(\tau)=\mathbb{P}({\cal N}(\tau)=n)$.
Then, by conditioning arguments we can construct from (\ref{fE}) the discrete-time density for the first arrival $E_t =J-\tau$ in the DTARP ${\cal N}_{\tau}(t)$ by summing up over all $n$ to arrive at
\beq
\label{conditional_series}
f_E(\tau,t) = \sum_{n=0}^{\infty} f_E(\tau,t,n).
\eeq
This summation stops at 
$n=\tau$ (as $f_E(\tau,t,n)=0$ for $n>\tau$) and
with the (double-) generating function
\beq
\label{double_gen_density}
 {\bar f}_E(w,u) =  \sum_{n=0}^{\infty}  {\bar f}_E(w,u,n) = 
 \frac{u}{u-w}\frac{{\bar \psi}(u)-{\bar \psi}(w)}{1-{\bar \psi}(w)} \hspace{1cm} (|u|\leq 1,\hspace{0.5cm} |w|<1,\hspace{0.5cm}  u\neq w).
\eeq
It is worthy of mention that for $w=0$ this recovers the GF of the waiting time density of the original process
${\bar f}_E(0,u)={\bar \psi}(u)$ as well as ${\cal N}_{\tau=0}(t)={\cal N}(t)$ recovers the original counting process. 
For later use we need to consider the case $w=u$ of this GF which is defined by 
$$
\lim_{w\to u}\frac{{\bar \psi}(u)-{\bar \psi}(w)}{u-w} = \frac{d {\bar \psi}(u)}{du}
$$
thus 
\beq\label{w_equal_u}
{\bar f}_E(u,u) = \frac{u \frac{d {\bar \psi}(u)}{du}}{1-{\bar \psi}(u)}
\eeq
and in the same way all cases $w=u$ are subsequently defined.
We notice that $u \frac{d {\bar \psi}(u)}{du}$ is the GF of $t\psi(t)$.
In view of the (absolute) monotonicity of ${\bar \psi}(u) > {\bar \psi}(w) $ for $u>w$
we confirm that (\ref{double_gen_density}) and (\ref{w_equal_u}) are for $w,u \in [0,1)$ non-negative.
Putting $u=1$ yields ${\bar f}_E(w,1) = \frac{1}{1-w}$ reflecting the normalization 
$\sum_{t=0}^{\infty}f_E(\tau,t) = 1$, i.e. $f_E(\tau,t)$ is a density on $t$ [Remark:  Eq. (\ref{double_gen_density}) is the discrete version of Eq. (6.2) in \cite{GodrecheLuck2001}].
\\[1mm]
Consider for a moment the (memoryless) Bernoulli process ${\cal N}_B(t)$ with ${\bar \psi}_{Ber}(u)= pu/(1-qu)$. This yields for (\ref{double_gen_density})
\beq
\label{bernoulli}
{\bar f}_E(w,u)_{Ber} = \frac{{\bar \psi}(u)_{B} }{1-w} =\frac{1}{1-w} \frac{pu}{1-qu} 
\eeq
corresponding to the unchanged geometric Bernoulli waiting-time PDF $f(\tau,t)_{Ber}= pq^{t-1}$, independent of the aging parameter $\tau$ due to
the Markovian nature of the Bernoulli process.
\\[1mm]
To demonstrate the aging effect in a non-Markovian renewal process it appears instructive to consider (\ref{double_gen_density}) for Sibuya waiting times with ${\bar \psi}_{\mu}(u)=1-(1-u)^{\mu}$ ($\mu \in (0,1)$) and $u,w \to 1$
to see the asymptotics for large $t$ and $\tau$ using Tauberian arguments.
Letting first $w\to 1-$ leads with (\ref{double_gen_density}) to the asymptotics $${\bar f}_{E}(u,w) \sim u\frac{1-{\bar \psi}(u)}{1-u} \frac{1}{1-{\bar \psi}(w)}  \sim (1-u)^{\mu-1}(1-w)^{-\mu}  , \hspace{1cm} (\, 1-w \ll 1-u \to 0\,). $$
This yields, as leading contribution in $\tau$, the scaling
\beq
\label{large-time-asymp-Sib}
f_{E,\mu}(t,\tau) \sim 
\frac{\Gamma(t+1-\mu)}{\Gamma(t+1)\Gamma(1-\mu)} \frac{\Gamma(\tau+\mu)}{\Gamma(\tau+1)\Gamma(\mu)}    \sim \frac{1}{\Gamma(1-\mu)\Gamma(\mu)} \frac{\tau^{\mu-1}}{t^{\mu}} , \hspace{1cm} (\tau \gg t \to\infty)
\eeq
which is the large time limit 
for strong aging $\tau \gg t \gg 1$ in the density of the first event in the (Sibuya-) DTARP. This type of scaling holds in general when ${\bar \psi}(t)$ is fat-tailed.
The aging effect decreases with a $\tau^{\mu-1}$ power-law with the aging parameter. This relation is in agreement with the result reported by Barkai for continuous times (\cite{Barkai2003}, Eq. (9)) for the strong aging (highly aged) limit $t_a\gg t_1$ (in his notation and identify $t_a=\tau$ and $t_1=t$)
 and use Euler's reflection formula $\Gamma(\mu)\Gamma(1-\mu) = \frac{\pi}{\sin(\pi\mu)}$. We refer also to the discussion in
 \cite{SchulzBarkaiMetzler2014}.
\\[3mm]
Let us continue to consider the general case.
To obtain the state probabilities that $\mathbb{P}[{\cal N}_{\tau}(t)=m]$, we have to take into account that after the first event $N_{\tau}(t)=1$ (governed by above density $f_E(\tau,t)$) the process is further developing with the IID waiting times of the original counting process.
Therefore, 
\beq
\label{genera_fuaging}
\mathbb{P}[{\cal N}_{\tau}(t) = m] =\Phi_{\tau}^{(m)}(t)= \sum_{r=0}^t f_E(\tau,r)\,\Phi^{(m-1)}(t-r) , \hspace{1cm} m=1,2,\ldots
\eeq
where $\Phi^{(r)}(k) = \mathbb{P}[{\cal N}(k)=r)]$ are the state probabilities of the original counting process (\ref{discrete-time}). We observe the initial condition $\Phi_{\tau}^{(m)}(0)=\delta_{m,0}$ reflecting 
${\cal N}_{\tau}(0)={\cal N}(\tau)-{\cal N}(\tau)=0$.
Formula (\ref{genera_fuaging}) has the double generating function
\beq
\label{double}
{\bar \Phi}_{w}^{(m)}(u) = {\bar f}_E(w,u) {\bar \Phi}^{m-1}(u) = 
{\bar f}_E(w,u) [{\bar \psi}(u)]^{m-1} \frac{1-{\bar \psi}(u)}{1-u} , \hspace{1cm} m=1,2,\ldots
\eeq
For $m=0$ we have the inversion of the double generating function
\beq
\label{double-gen-survival}
{\bar \Phi}_{w}^{(0)}(u) = \frac{1}{(1-u)(1-w)}- \frac{{\bar f}_E(w,u)}{1-u} , \hspace{1cm} m=0
\eeq
yields the survival probability 
$\mathbb{P}[{\cal N}_{\tau}(t)=0]$ 
[Remark: The last two Eqs. are the discrete-time counterparts of
Eqs. (3), (4) in \cite{SchulzBarkaiMetzler2014}].
\\[3mm]
Relations (\ref{double}) and (\ref{double-gen-survival}) can also be derived by the following
considerations.
Consider the conditional probability that ${\cal N}_{\tau}(t)=m$ given ${\cal N}(\tau)=n$:
\beq
\label{con_probability_m_nonzero}
\mathbb{P}[{\cal N}_{\tau}(t)=m|{\cal N}(\tau)=n]
=\Phi^{m,n}_{\tau}(t) =\langle \Theta(J_n\tau,J_{n+1})\Theta(J_{n+m}-\tau,t,J_{n+m+1}-\tau)\rangle  , \hspace{2mm} m=1,2,\ldots.
\eeq
Then consider separately the case ${\cal N}_{\tau}(t)=0$ when ${\cal N}(\tau)=n$ which writes
\beq
\label{condprob_mzero}
\mathbb{P}[{\cal N}_{\tau}(t)=0|{\cal N}(\tau)=n]
=\Phi^{0,n}_{\tau}(t) = \langle \Theta(J_n\tau,J_{n+1})\Theta(J_{n+1}-1-\tau-t) \rangle 
\eeq
where $\Theta(J_{n+1}-1-\tau-t)=1$ when $J_{n+1}>t+\tau$.
Then in the same way as above we compute the double generating functions
of the last two probabilities (nota bene $\tau \leftrightarrow w$ and $t \leftrightarrow u$) to arrive at
\beq
\label{con_probability_m_nonzero_genfun}
{\bar \Phi}^{m,n}_w(u)= {\bar \psi}(w)^n\frac{u}{u-w}
[{\bar \psi}(u)-{\bar \psi}(w)]{\bar \psi}(u)^{m-1}\frac{1-{\bar \psi}(u)}{1-u} ,  \hspace{1cm} m=1,2,\ldots
\eeq
and for $m=0$ we have
\beq
\label{con_probability_m_zero_genfun}
{\bar \Phi}^{0,n}_w(u) = \frac{{\bar \psi}(w)^n}{(1-u)}
\left[\frac{1-{\bar \psi}(w)}{1-w}
-\frac{u}{u-w}({\bar \psi}(u)-{\bar \psi}(w))\right]  
\eeq
where the summation over $n$ of ${\bar \Phi}^{m,n}_w(u)$
indeed yields double GFs (\ref{double}) and (\ref{double-gen-survival}), respectively, of the state probabilities $\mathbb{P}[{\cal N}_{\tau}(t)=m]$ of the DTARP. This reflects the conditional relation 
\beq
\label{conditioanl_relation}
\mathbb{P}[{\cal N}_{\tau}(t)=m] = 
\sum_{n=0}^{\infty}\mathbb{P}[{\cal N}_{\tau}(t)=m|{\cal N}(\tau)=n]
\eeq
where this series breaks at $n=\tau$ (as ${\cal N}(\tau) \leq \tau$). With the above GFs one can straight-forwardly also show that $\mathbb{P}[{\cal N}(\tau)=n] = \sum_{m=0}^{t}\mathbb{P}[{\cal N}_{\tau}(t)=m|{\cal N}(\tau)=n]$, recovering the state probabilities of the original counting process.
\subsection{Some pertinent DTARP generating functions}
\label{pertinent_DTARP}
Of interest for the evaluation of the mean square displacement (\ref{msd_discrte-time}) are averages of the type
\beq
\label{average_aging}
g_v(\tau,t) = \left\langle v^{{\cal N}_{\tau}(t)} \right\rangle = \sum_{m=0}^{\infty}
\mathbb{P}[{\cal N}_{\tau}(t)=m] v^m , \hspace{1cm} |v| \leq 1
\eeq
for a DTARP (\ref{agin_renewal}). For
$v=1$ (\ref{average_aging}) yields unity (normalization)
$$g_1(\tau,t)=\sum_{m=0}^{\infty}\mathbb{P}[{\cal N}_{\tau}(t)=m]=1.$$ These series stop at $m=t$ 
with $\mathbb{P}[{\cal N}_{\tau}(t)=m]=0$ for $m>t$ since ${\cal N}_{\tau}(t) \leq t$ 
(reflected by the feature ${\bar \Phi}_w^{(m)}(u)=O(u^m)$).
Therefore, $g_v(\tau,t)$ is a polynomial of degree $t$ in $v$ (`state polynomial of the DTARP').
We also mention that
\beq
\label{auto-correlation}
g_{-1}(\tau,t) = \langle (-1)^{{\cal N}(t+\tau)-{\cal N}(\tau)} \rangle =
 \left\langle (-1)^{{\cal N}(t+\tau)+{\cal N}(\tau)} \right\rangle = 
 \left\langle \sigma_{t+\tau}\sigma_{\tau} \right\rangle , \hspace{1cm} t, \tau >0
\eeq
is the auto-correlation function of the steps for non-zero $t,\tau$.
We will come back to these properties by means of generating functions. Noteworthy is that $g_v(\tau,t)\big|_{\tau=0} =  \langle v^{{\cal N}_0(t)} \rangle = {\cal P}(v,t)$ recovers the state polynomial of the original counting process 
${\cal N}_0(t)={\cal N}(t)$.
We are interested in the triple generating function
\beq
\label{double_gen_aging}
\begin{array}{cllr} 
\ds {\bar g}_v(w,u) & = \ds \sum_{\tau=0}^{\infty}\sum_{t=0}^{\infty} w^{\tau}u^t  \langle v^{{\cal N}_{\tau}(t)} \rangle , &|w|,|u| < 1 , & |v| \leq 1\\ \\ & = \ds \sum_{m=0}^{\infty} \sum_{\tau=0}^{\infty}\sum_{t=0}^{\infty}v^m w^{\tau}u^t \mathbb{P}[{\cal N}_{\tau}(t)=m] , & u\neq w &\\ \\
 & = \ds \sum_{m=0}^{\infty}{\bar \Phi}_w^{(m)}(u)v^m  & & \\ \\ 
 & = \ds \frac{1}{(1-w)(1-u)} - \frac{(1-v)u}{(1-u)(u-w)[1-v{\bar \psi}(u)]}\frac{[{\bar \psi}(u)-{\bar \psi}(w)]}{[1-{\bar \psi}(w)]} & & \\ \\
 & = \ds   \frac{1}{(1-w)(1-u)} - \frac{(1-v){\bar f}_{E}(w,u)}{(1-u)[1-v{\bar \psi}(u)]} & &
 \end{array}
\eeq
where we used (\ref{double}) and (\ref{double-gen-survival}) with (\ref{double_gen_density}). We directly confirm
${\bar g}_v(w,u)\big|_{v=1} = \frac{1}{(1-w)(1-u)}$ reflecting $g_1(\tau,t)=1$.
For $v=-1$, formula (\ref{double_gen_aging}) is the GF
of the expected value $g_{-1}(\tau,t) = \langle (-1)^{{\cal N}_{\tau}(t)} \rangle$ for
the DTARP. Note that formula (\ref{double_gen_aging}) for $w=0$ yields the GF of the state polynomial of the original process (as ${\cal N}_{\tau=0}(t)={\cal N}(t)$):
\beq\label{wnull}
{\bar g}_v(0,u) = {\bar {\cal P}}(v,u) = \sum_{t=0}^{\infty}u^t  \langle v^{{\cal N}(t)} \rangle
= \frac{1-{\bar \psi}(u)}{1-u} \frac{1}{1-v{\bar \psi}(u)}.
\eeq
For the Bernoulli process we get for (\ref{double_gen_aging}),
\beq
\label{Bernoulli_v}
{\bar g}_{v}(w,u)_{Ber} = 
\frac{1}{(1-w)\{1-u(1+p[v-1])\}}
\eeq
thus
\beq
\label{Bernoulli_no_aging}
{\bar g}_{-1}(w,u)_{Ber} =  \frac{1}{(1-w)[1-u(1-2p)]}
\eeq
and hence $g_{-1}(\tau,t)_{Ber}=(1-2p)^t =\langle (-1)^{{\cal N}_B(t)}\rangle$ which is our result of Section \ref{Bernoulli}
and is independent of the aging parameter $\tau$,
reflecting the Markovian nature of the Bernoulli process.
\\[3mm]
In order to evaluate the mean square displacement (\ref{msd_discrte-time}) we need to consider the double GFs of functions of the form
$h(\tau,t)=\sum_{k=0}^{t-\tau}g(\tau,k)$ ($t\geq \tau$), namely
\beq
\label{doublegen}
\begin{array}{clr}
{\bar h}(w,u) & = \ds  \sum_{t=0}^{\infty}u^t\sum_{\tau=0}^{\infty}w^{\tau} \Theta(t-\tau) \sum_{k=0}^{t-\tau} g(\tau,k)  = \ds  \sum_{t=0}^{\infty} u^t H_w(t) \\ \\
 & = \ds  \sum_{\tau=0}^{\infty} (uw)^{\tau}\sum_{s=0}^{\infty} u^s  \sum_{k=0}^{s}  g(\tau,k) & \\ \\
 & = \ds \sum_{s=0}^{\infty} u^s  
 \sum_{k=0}^{s} \Theta(s-k)\left[\sum_{\tau=0}^{\infty} (uw)^{\tau} g(\tau,k)\right] & \\ \\
 & = \ds \frac{{\bar g}(uw, u)}{1-u} = {\bar H}_w(u)
\end{array}
\eeq
where we have performed summation over the new  indices $s=t-\tau$ and $\tau$. Be reminded of our notation for double GFs ${\bar g}(a,b) =\sum_{r=0}^{\infty}\sum_{s=0}^{\infty} a^r b^s g(r,s)$.
We introduced in the first line of this relation the auxiliary quantity
\beq
\label{auxiliary_quant}
H_w(t) = \sum_{\tau=0}^{\infty}w^{\tau} \Theta(t-\tau) \sum_{k=0}^{t-\tau} g(\tau,k)
\eeq
(breaking at $\tau=t$).
Next, we focus on $g(\tau,t)=g_v(\tau,t)= 
\langle v^{{\cal N}_{\tau}(t)} \rangle $.
Using (\ref{average_aging}), (\ref{double_gen_aging}) with (\ref{doublegen}) we get for the double GF of
$h_v(\tau,t) = \sum_{k=0}^{t-\tau}\langle v^{{\cal N}_{\tau}(k)}\rangle$ the expression
\beq
\label{double-gener}
\begin{array}{clr}
\ds {\bar h}_v(w,u) & =  \ds \sum_{\tau=0}^{\infty}\sum_{t=0}^{\infty}w^{\tau}u^t \Theta(t-\tau) \sum_{k=0}^{t-\tau}\langle v^{{\cal N}_{\tau}(k)}\rangle  = \frac{{\bar g}_v(uw, u)}{1-u} & \\ \\ 
 & = \ds  \frac{1}{(1-wu)(1-u)^2} - \frac{(1-v){\bar f}_{E}(wu,u)}{(1-u)^2[1-v{\bar \psi}(wu)]} & \\ \\ 
 & = \ds  \frac{1}{(1-wu)(1-u)^2} - \frac{(1-v)[{\bar \psi}(u)-{\bar \psi}(uw)]}{(1-w)(1-u)^2[1-v{\bar \psi}(u)][1-{\bar \psi}(uw)]}.
 \end{array}
\eeq
For the following evaluation it is useful to take into account that ${\bar f}(u)/(1-uw)$ is the double GF of $\Theta(t-\tau) f(t-\tau)$. 
This double GF contains for $w=1$ the GF of $\sum_{k=0}^tf(k)$. 
Therefore,
\beq
\label{estab_GF}
\begin{array}{clr}
\ds {\bar h}_{-1}(1,u) &=
\ds \sum_{t=0}^{\infty} u^t \sum_{\tau=0}^{t}
\sum_{k=0}^{t-\tau} 
\langle (-1)^{{\cal N}_{\tau}(k)}\rangle &
\\ \\
 & = \ds \frac{1}{(1-u)^3} - \frac{2u \frac{d{\bar \psi}(u)}{du}}{(1-u)^2(1-[{\bar \psi}(u)]^2)}
\end{array}
\eeq
where 
in order to evaluate (\ref{double-gener}) for $w=1$ we accounted for (\ref{w_equal_u}), i.e.
$$
{\bar \psi}(uw)= {\bar \psi}[u- u(1-w)] =
{\bar \psi}(u) - u(1-w) \frac{d{\bar \psi}(u)}{du} + o(1-w) $$
thus
\beq
\label{hv_res_gq}
{\bar h}_v(1,u) = \frac{1}{(1-u)^3}-
\frac{(1-v)u}{(1-u)^2[1-v{\bar \psi}(u)][1-{\bar \psi}(u)]}\frac{d {\bar \psi}(u)}{du}
\eeq
and ${\bar h}_v(1,0)=1$.
\\[3mm]
For the Bernoulli process we have
\beq
\label{Bernoulli_hs}
\begin{array}{clr}
\ds {\bar h}_{-1}(1,u)_{Ber} & = \ds \frac{1}{(1-u)^3}-
\frac{2pu}{(1-u)^3[1-u(1-2p)]} & \\ \\
 & = \ds \frac{1}{(1-u)^2[1-u(1-2p)]} &
 \end{array}
\eeq
consistent with our previous results, namely
\beq
\label{previous}
\begin{array}{clr} 
\ds {\bar h}_{-1}(1,u)_{Ber} & =\ds \sum_{t=0}^{\infty}u^t \sum_{\tau=0}^t\sum_{k=0}^{t-\tau}\langle (-1)^{{\cal N}_B(k)}\rangle &  \\ \\
& = \ds \sum_{t=0}^{\infty}u^t \sum_{\tau=0}^t\sum_{k=0}^{\tau}\langle (-1)^{{\cal N}_B(k)}\rangle &
\end{array}
\eeq
where
$$\sum_{t=0}^{\infty}u^t\langle (-1)^{{\cal N}_B(t)}\rangle = \sum_{t=0}^{\infty}u^t(1-2p)^t = \frac{1}{1-u(1-2p)} 
$$
reflecting again the Markovian nature of Bernoulli.
\\[3mm]
In order to evaluate the mean square displacement (\ref{msd_discrte-time}) of the SRW we
need to remove the term $\tau=0$ in (\ref{estab_GF}) and define its GF (summation starting at $\tau=1$ and, considering (\ref{wnull}), we take into account that ${\cal N}_0(t)={\cal N}(t)$)
\beq
\label{genfu}
\begin{array}{clr}
\ds {\bar K}(u) & = \ds  \sum_{t=1}^{\infty}u^t\sum_{\tau=1}^t \sum_{k=0}^{t-\tau}\langle (-1)^{{\cal N}_{\tau}(k)} \rangle &  \\ \\
& = \ds {\bar h}_{-1}(1,u) - \sum_{t=0}^{\infty}u^t \sum_{k=0}^{t} \langle (-1)^{{\cal N}(k)}\rangle = {\bar h}_{-1}(1,u) -{\bar h}_{-1}(0,u) & \\ \\
 & = \ds  \frac{1}{(1-u)^3} - \frac{2u \frac{d{\bar \psi}(u)}{du}}{(1-u)^2(1-[{\bar \psi}(u)]^2)} - \frac{1-{\bar \psi}(u)}{(1-u)^2} \frac{1}{1+{\bar \psi}(u)}        & \\ \\
 & = \ds \frac{1}{(1-u)^3} -\frac{1}{(1-u)^2(1-[{\bar \psi}(u)]^2)}\left(2u\frac{d{\bar \psi}(u)}{du}+[1-{\bar \psi}(u)]^2\right) &
\end{array}
\eeq
where the necessary property 
${\bar K}(u)\big|_{u=0}=K(t)\big|_{t=0}=0$ 
is fulfilled. Now, to obtain the GF (which we denote by ${\bar X}^{(2)}(u) = \sum_{t=0}^{\infty}u^t \langle\,  X_t^2\, \rangle$) of the mean square displacement
we have to take into account that 
$$\langle \, X_t^2 \, \rangle =2K(t)-t$$ 
where $t$ has GF $\frac{u}{(1-u)^2}$,
thus
\beq
\label{GF_msd_general}
{\bar X}^{(2)}(u) = 2{\bar K}(u)-\frac{u}{(1-u)^2}.
\eeq
\\[3mm]
For Bernoulli this yields 
\beq
\label{Bernoulli_K}
{\bar K}_B(u)   
= \frac{u}{(1-u)^2[1-u(1-2p)]} =  u {\bar h}_{-1}(1,u),
\eeq
thus the GF of the mean square displacement (\ref{variance_bernoulli}) yields
\beq
\label{Ber_msdgenfu}
{\bar X}^{(2)}(u)_{Ber} = 2{\bar K}_B(u)
-\frac{u}{(1-u)^2} = \frac{u}{(1-u)^2}\frac{[1+u(1-2p)]}{[1-u(1-2p)]},
\eeq
in agreement with the results derived in Section \ref{Bernoulli}. For the unbiased case 
$p=\frac{1}{2}$ we have ${\bar X}^{(2)}(u)_{Ber}= \frac{u}{(1-u)^2}$ and
with 
$(\langle X_t\rangle_B))_{p=\frac{1}{2}} =0$ (see Eq. (\ref{expected_pos_Bern})) thus we recover
$\langle X^2_{p=\frac{1}{2}}\rangle_B= {\cal V}_{p=\frac{1}{2}}(t)=t$. 
Let us check the deterministic limiting case $p=0+$,
which corresponds to the strictly increasing walk with
$\langle X_t^2\rangle =t^2$.
Formula (\ref{Ber_msdgenfu}) then becomes
\beq
\label{p_equal_zero}
{\bar X}^{(2)}(u)_{Ber} = \frac{u(1+u)}{(1-u)^3} ,\hspace{1cm} (p=0+).
\eeq
Then, account for the GFs 
$\frac{1}{1-u} \to \Theta(t)$, $\frac{1}{(1-u)^2} \to  t+1 $, 
$\frac{1}{(1-u)^3} \to \sum_{k=0}^t(k+1) =(t+2)(t+1)/2$, $u/(1-u)^3 \to t(t+1)/2$
and $u^2/(1-u)^3 \to (t-1)t/2$. With the last two GFs we recover the result
of Eq. (\ref{expand_p})
\beq
\label{resp_null}
\langle X^2_{p=0+}\rangle_B= \frac{1}{2} \left[t(t+1)+t(t-1)\right] = t^2
\eeq
and hence we have, with (\ref{delopital}), 
\beq
\label{Var_Ber_p0}
{\cal V}_{p=0+}(t)=0
\eeq
corresponding to the deterministic walk with (a.s.) unit steps
in ${\tilde \sigma}_0$-direction (see (\ref{delopital})).
\\[1mm]
Finally, for $p=1$, formula (\ref{Ber_msdgenfu}) yields 
${\bar X}^{(2)}(u)_{Ber}=\frac{u}{(1-u)(1+u)}$, corresponding to the deterministic 
oscillatory motion where the mean square displacement oscillates between zero and one with ${\bar X}^{(2)}(u)_{Ber} \sim \frac{1}{2(1-u)}$ as $u \to 1-$ and yields the
large time limit $\langle X_t^2\rangle_B \sim \frac{1}{2}$ (see Section \ref{Bernoulli}).
\subsection{Causal distributions and their Laplace transforms}
\label{Lapl_Trafo}
We deal with causal functions and distributions $\Theta(t)f(t)$ ($t\in \mathbb{R}^{+}$ and $\Theta(t)$ indicates the Heaviside step function with $\Theta(t)=1$ for $t\geq 0$ and $\Theta(t)=0$ elsewhere). We introduce the Laplace transform as
\beq
\label{Laplace_trafo}
{\hat f}(s)= \int_{0-}^{\infty}e^{-st}\Theta(t)f(t){\rm d}t
\eeq
with suitably chosen Laplace variable $s$. We mention that
$s^m{\hat f}(s)$ is the Laplace transform of $\frac{d^m}{dt^m}[\Theta(t)f(t)]$ as all boundary terms
are vanishing at $0-$. As a consequence $\Theta(0-)=0$.
\subsection{Discrete-time Prabhakar kernel}
\label{Discrete_Prabhakar}
We have represented
some of our results in terms of discrete-time Prabhakar kernels.
The name `Prabhakar kernel' comes from the fact that these kernels involve
`Prabhakar functions' first introduced in \cite{Prabhakar1971}
as generalization of the Mittag-Leffler function. Consult \cite{Giusti-et-al-Prabhakar2020,Giusti2020} for a general outline of Prabhakar fractional calculus and pertinent applications.
We recall here some of the essential features of their discrete-time versions, the discrete-time Prabhakar kernels. We introduced the discrete-time versions recently \cite{MichelitschPolitoRiascos2021,TMM_FP_APR_Prab_Warsaw2021}.
The discrete-time Prabkakar kernel which we denote with
$p_{\mu,\nu}^{\,\gamma}(\lambda,t)$ is defined by its GF (\ref{PrabGF})
\beq
\label{PrabGF_pp}
{\bar p}_{\mu,\nu}^{\,\gamma}(\lambda,u) =
\frac{(1-u)^{-\nu}}{(1-\lambda(1-u)^{-\mu})^{\gamma}} , \hspace{1cm} |u| < 1.
\eeq
Here it is sufficient to consider 
the range $\mu \in (0,1]$ 
with $\nu>0$, $\lambda >0$ and $\gamma=1$ where we write
${\bar p}_{\mu,\nu}^1(\lambda,u) ={\bar p}_{\mu,\nu}(\lambda,u)$. Then we have
\beq
\label{expand_small_lam}
{\bar p}_{\mu,\nu}(\lambda,u) = \left\{\begin{array}{cr} \ds  \sum_{m=0}^{\infty}\lambda^m(1-u)^{-(\nu+m\mu)}, & \hspace{0.5cm} \ds  \lambda|1-u|^{-\mu} < 1\\ \\ 
 \ds - \sum_{m=1}^{\infty}\lambda^{-m}(1-u)^{\mu m-\nu}, & \hspace{0.5cm} \ds \lambda|1-u|^{-\mu} > 1 \end{array} \right.
\eeq
(where $\lambda(1-u)^{-\mu} \neq 1$).
Thus we get for the Prabhakar kernel the expansion (where we write $p_{\mu,\nu}(\lambda,t) = p_{\mu,\nu}^1(\lambda,t)$) 
\beq
\label{Prabhakar_funct}
p_{\mu,\nu}(\lambda,t) = \frac{1}{t!}\frac{d^t}{du^t}{\bar p}_{\mu,\nu}(\lambda,u)\big|_{u=0} = \left\{\begin{array}{lr} \ds  \frac{1}{t!}\sum_{m=0}^{\infty} 
\lambda^m \frac{\Gamma(\nu+m\mu+t)}{\Gamma(\nu+m\mu)} , &\ds \hspace{0.5cm}   |\lambda| < 1 \\ \\ \ds
\frac{(-1)^{t+1}}{t!} \sum_{m=1}^{\infty} \lambda^{-m}\frac{\Gamma(\mu m-\nu+1)}{\Gamma(\mu m-\nu+1-t)} , & \ds |\lambda| > 1 \end{array}\right.\hspace{-1cm} t \in \mathbb{N}_0
\eeq
which converges absolutely for finite $t$ and is divergent for $\lambda=1$. Note that the coefficients (Pochhammer symbol) $(c)_k=\Gamma(c+k)/\Gamma(c)
=c(c+1) \ldots (c+k-1)$ ($k\in \mathbb{N}_0$)  fulfill $(0)_k=0$ for $k\geq 1$) and for $k=0$ we have $(c)_0=1$ and especially $(0)_0=1$.
For $t=0$ we verify
$p_{\mu,\nu}(\lambda,0) ={\bar p}_{\mu,\nu}(\lambda,0) = 1/(1-\lambda)$. For the exploration of the case $\lambda>1$ (relevant in our evaluation where $\lambda=2$), it may be convenient
to implement the product representation
\beq
\label{auxiliary}
\begin{array}{clr}
\ds (-1)^t\frac{\Gamma(\mu m-\nu+1)}{\Gamma(\mu m-\nu+1-t)} & =(\nu-m \mu)_t = 
\frac{\Gamma(\nu-m \mu+t)}{\Gamma(\nu-m\mu)}.
\end{array}
\eeq
Evoking Tauberian arguments and accounting for
$\Gamma(c+t)/\Gamma(t+1) \sim t^{c-1}$ ($t\to \infty$) we have
$$
\frac{(-1)^{t+1}}{t!}\lambda^{-m}\frac{\Gamma(\mu m-\nu+1)}{\Gamma(\mu m-\nu+1-t)} = -\lambda^{-m}\frac{\Gamma(t+\nu-m\mu)}{\Gamma(\nu-m\mu)\Gamma(t+1)} \sim -\lambda^{-m} \frac{t^{\nu-m\mu-1}}{\Gamma(\nu-m\mu)} , \hspace{0.5cm} (t \to \infty)
$$
where the term for $m=1$ is the dominant. We can extract it from (\ref{expand_small_lam}) (with $u\to 1-$) and is identical with the limit $\lambda \gg 1$:
\beq
\label{as_t_small_and_large}
p_{\mu,\nu}(\lambda,t) 
 \sim -\frac{1}{\lambda} \frac{t^{\nu-\mu-1}}{\Gamma(\nu-\mu)}  ,\hspace{1cm} (\nu \neq \mu, \hspace{0.5cm} t\to \infty).
\eeq
For $\nu=\mu$ the term $m=1$ yields $-\lambda^{-1}\delta_{t0}$, thus the tail is dominated by the next term ($m=2$), namely
\beq
\label{mu_mu}
p_{\mu,\mu}(\lambda,t) \sim -\frac{1}{\lambda^2} \frac{\Gamma(t-\mu)}{\Gamma(-\mu)\Gamma(t+1)} \sim  -\frac{1}{\lambda^2}\frac{t^{-\mu-1}}{\Gamma(-\mu)} , \hspace{1cm} (t \to \infty).
\eeq
Be reminded that for $r \in \mathbb{N}$, expressions of the form $\frac{t^{-r-1}}{\Gamma(-r)}$ are defined 
in the Gel'fand-Shilov sense \cite{Gelfand-Shilov1968} as $\frac{t^{-r-1}}{\Gamma(-r)} = 
\frac{d^r}{dt^r}\delta(t) =0$ (being strictly null for large times).
Thus, in (\ref{Kmu_exact}) the Prabhakar kernels then have the large time asymptotics: 
\beq
\label{asymp_Prabhakar}
p_{\mu,\mu+2}(2,t) \sim -\frac{t}{2}, \hspace{0.1cm}
p_{\mu,\mu+3}(2,t)\sim -\frac{t^2}{4},\hspace{0.1cm} p_{\mu,2}(2,t)\sim -\frac{t^{1-\mu}}{2\Gamma(2-\mu)}, \hspace{1cm}  (t\to\infty). 
\eeq
It is worthy of mention that the following scaling limit exists \cite{MichelitschPolitoRiascos2021,TMM_FP_APR_Prab_Warsaw2021}
\beq
\label{scaling_lim}
e^{1}_{\mu,\nu}(\lambda_0,t) = \lim_{h\to 0} h^{\nu-1} 
p_{\mu,\nu}(\lambda_0h^{\mu},\frac{t}{h}) 
= t^{\nu-1}\sum_{m=0}^{\infty} \frac{(\lambda_0t^{\mu})^m}{\Gamma(\nu+m\mu)} = t^{\nu-1}E_{\mu,\nu}(\lambda_0t^{\mu})
, \hspace{1cm} (t \in h{\mathbb N} \to \mathbb{R}^+)
\eeq
where $e^{\gamma=1}_{\mu,\nu}(\lambda_0,t)$ is the continuous-time
Prabhakar kernel
\cite{Giusti-et-al-Prabhakar2020,Giusti2020}.
In (\ref{scaling_lim}) we use the expansion (\ref{Prabhakar_funct}) for $\lambda=\lambda_0h^{\mu} < 1$. This scaling limit contains the two parameter Mittag-Leffler function (\ref{two-param_ML})
and a new parameter $\lambda_0>0$ (of physical dimension $[time^{-\mu}]$ and independent of $h$).
Therefore, indeed (\ref{Prabhakar_funct}) is a discrete-time approximation of
the (Prabhakar-) kernel
(\ref{scaling_lim}).
The Laplace transform of the continuous-time Prabakar kernel (\ref{scaling_lim}) is connected with
(\ref{PrabGF}) by the scaling limit $\lim_{h\to 0} h^{\nu}{\bar p}_{\mu,\nu}(\lambda_0h^{\mu},e^{-hs}) = 
s^{-\nu}/(1-\lambda_0s^{-\mu})= {\hat e}_{\mu,\nu}^1(\lambda_0,s)$. 
\end{appendix}
\subsection*{Acknowledgments}
F.~Polito has been partially supported by INdAM/GNAMPA.

\end{document}